\documentclass[review,onefignum,onetabnum]{article}

\usepackage[utf8]{inputenc}
\usepackage[scaled]{helvet} % Helvetica.  Similar to Arial as specified by NERC
\usepackage{amsmath,amssymb,amsthm}
\usepackage{subcaption}
\usepackage{amsfonts}
\usepackage{mdframed}
\usepackage{dsfont}
\usepackage{tikz}
\usepackage{nicefrac}
\usepackage[margin=2.5cm]{geometry}
\usepackage[many]{tcolorbox}
\usepackage{array}
\usepackage{empheq}
\usepackage{siunitx}
\usepackage{placeins}
\newtheorem{theorem}{Theorem}
\newtheorem{lemma}{Lemma}
\newtheorem*{remark}{Remark}

\usepackage{graphicx} % Includes graphics
\usepackage{float} % Allows forcing of figure locations
\usepackage{wrapfig} % Wraps text round images
\usepackage[font=small,format=plain,labelfont=bf,textfont=it]{caption}
\usepackage{subcaption}
\usepackage{multicol}
\usepackage{enumerate}
\usepackage{hyperref}
\usepackage{authblk}
\usepackage[capitalise]{cleveref}

\newcommand{\R}{\mathbb{R}}

\newcommand{\di}[1]{\,\mathrm{d}#1}
\newcommand{\ID}{\operatorname{Id}}

\makeatletter
\newcommand{\tline}{%
    \noalign {\ifnum 0=`}\fi \hrule height 1pt
    \futurelet \reserved@a \@xhline
}
\newcolumntype{"}{@{\hskip\tabcolsep\vrule width 1pt\hskip\tabcolsep}}
\makeatother
\title{Analysis of a Hydrodynamic Thermoelastic Mixture Model}

\author[$\dagger$]{Michael Eden}
\author[$\star$]{Meraj Alam}
\author[$\ddag$]{Prakash Kumar}
\author[$\S$]{Raja GP Sekhar}

\affil[$\dagger$]{University of Regensburg, Regensburg, Germany}
\affil[$\star$]{Mahindra University, Hyderabad, India}
\affil[$\ddag$]{SRM University AP, Amaravati, India}
\affil[$\S$]{IIT Kharagpur, Kharagpur, India}
\begin{document}

\maketitle
%%%%%%%%%%%%%%%%%%%%%%%%%%%%%%%%%%%%%%%%%%%%%%%%%%%%%%%%%%%%%%%%%
% Abstract
%%%%%%%%%%%%%%%%%%%%%%%%%%%%%%%%%%%%%%%%%%%%%%%%%%%%%%%%%%%%%%%%%

This study proposes and explores a hydrodynamic thermo-elasticity system within the context of mixture models, encompassing fluid and solid phases, with particular emphasis on biological tissues and tumor-related phenomena.
While tumor growth dynamics are not explicitly modeled yet, this work examines the interaction between thermal effects and hydrodynamics on short-time scales where the tumor size typically remains stable.
We establish the existence of a unique weak solution within the framework of implicit evolution equations where we exploit the intricate coupling mechanisms inherent to the PDE system.
To further investigate the model, we then study the one-dimensional model and explore in detail the complex interplay between fluid flow, solid deformation, and heat transfer.
This complex coupled system of equations is then reduced over a short time scale to obtain semi-analytical solutions.

{\bf Keywords: Mixture theory, thermoelasticity, modeling of biological tissue, coupled PDE systems}  

{\bf MSC2020: 74E30, 74F05, 74F10, 35M30}

%%%%%%%%%%%%%%%%%%%%%%%%%%%%%%%%%%%%%%%%%%%%%%%%%%%%%%%%%%%%%%%%%
% Introduction
%%%%%%%%%%%%%%%%%%%%%%%%%%%%%%%%%%%%%%%%%%%%%%%%%%%%%%%%%%%%%%%%%
\section{Introduction}

%\begin{itemize}
   % \item Extended literature review
   % \begin{itemize}
    %    \item \cite{Ezzat2020-ak}: Thermo-mechanical memory responses of biological viscoelastic
   %            tissue with variable thermal material properties
    %    \item \cite{Zhang2022-et}: Bio-thermo-viscoelastic behavior in %multilayer skin tissue
    %    \item \cite{Ciarletta2013-zq}: Mechano-transduction in tumour growth modelling
     %   \item \cite{Majaron1999-yc}: Thermo-mechanical laser ablation of soft biological tissue:
     %          modeling the micro-explosions
               %%%%%%%%%%%%%%%%%%%%%%%LiteratureAddedbyMeraj%%%%%%%%%%%%%%%%%%%%%%%%%
      %  \item \cite{he2010local}:  A LOCAL THERMAL NONEQUILIBRIUM %POROELASTICTHEORY FOR FLUID SATURATED POROUS MEDIA
      %  \item \cite{li2019thermomechanical} Thermomechanical response of porous biological tissue based on local thermal non-equilibrium
      %  \item       \cite{cieslak2022model}:                        A model in one-dimensional thermoelasticity
   % \end{itemize}
%\end{itemize}

Modeling the mechanical behavior of biological tumors using porous medium equations is very common
\cite{astanin2008multiphase,preziosi2009multiphase}.
However, considering the elastic nature of the tumor structure, there are two commonly known approaches for modeling soft tissue, namely the theory of poroelasticity (we point to \cite{coussy2004}) and the theory of mixture (e.g., \cite{B76}).
In the former, soft tissues are supposed to exhibit sponge-like behavior as fluid-saturated elastic porous solids, whereas in the latter, the continuous binary mixture of the solid (s) and fluid (f) phases represents the dynamics.
To estimate the temperature effect on solid deformations, thermal expansion effects are to be considered in the constitutive equations, and these equations are termed thermoporoelasticity theory \cite{Andreozzi2019-hh}.
Detailed literature review indicates that the scientific community has only recently begun showing enthusiasm to understand the effects of thermal expansion on fluid transport in tissues.

Pioneering work on mixtures as continua interacting can be assigned to Truesdell and Toupin \cite{truesdell1960theory} and Green and Naghdi \cite{GREEN1965231,GREEN1968631}. Consequently, several authors have considered the mathematical analysis involving mixtures of elastic and viscoelastic constituents.
Iesan \cite{iesan91} developed the theory of binary mixtures of thermoelastic solids where each of the phaes is treated as thermoelastic solid.
The corresponding dependent variables are displacement gradients, relative displacement, temperature and temperature gradient. 
They have shown a uniqueness theorem for the governing system of thermoelasticity for linear theory.
In \cite{Burchuladze00}, the authors considered diffusion model of the linear theory of binary mixtures of the thermoelastic solids.
They have considered a class of boundary value problems and shown uniqueness results.
Various composites that are used in several industrial applications occur as multiphase mixtures.
These can be treated as interacting continua.
The first thermomechanical theory for mixtures of Newtonian fluids in which the two constituents have different temperatures was presented by Eringen and Ingram \cite{eringen1965,INGRAM1967289}.
One can see the hierarchy of models in the context of mixture theory in \cite{RT95}.
\cite{FERNANDEZ20171886} considered the one dimensional problem of two thermoelastic mixtures with two different temperatures. 
The displacement of each of the phases is assumed to be governed by a linear elasticity equation (which is hyperbolic) with the relative displacement of the two phases appearing in individual phase equations obeying Newton's third law.
Further, the individual phase temperatures also appear in the solid phase displacement equation.
The temperature fields obey the energy balance equation which is parabolic in nature. Each of the energy balance equation is coupled to the temperature of the other phase and also the corresponding phase velocity of the solid phases.
The coupled system is reformulated as a variational problem and the corresponding existence and uniqueness results are discussed.
A finite element approximation of the problem is developed and some a priori error estimates are given followed by some numerical simulations.

Applying heat to human tissue is a well-known treatment mechanism for specific health conditions.
Experimental investigations show oscillatory temperature responses when supplying abrupt constant energy to tissue \cite{Yuan2008}.
The non-homogeneous nature of the tissue supports that the heat exchange between blood and tissues always happens at a finite speed.
This led to the development of different mathematical models to evaluate the thermal
behaviour taking place in living tissues \cite{Wang2021-uj}.
Estimating the thermoelastic deformation caused by the elevated temperature during thermal ablation procedures is very important to treat cancer. Typical stress–strain equations for an elastic material with thermal load can be thought of as an elementary model for this purpose, which is termed as thermoelastic wave equation \cite{Park2018-xo}. 
In \cite{Zhang2021}, nonlinear characteristics of bio-heat transfer in deformed soft tissues with thermal expansion/shrinkage are analyzed via a thermo-visco-hyperelastic model.
They developed an efficient solution method for solving the coupled model in different computational domains.
The solution method is then used to simulate thermal ablation in the liver
\cite{Wang2021-uj} have proposed a generalized thermo-elastic model involving a dual-phase-lag model of bioheat transfer to model the thermal response of skin tissue.
The thermo-elastic response of skin tissue subjected to sudden heating on its boundary is solved by this analytical approach.
They have evaluated the contribution of heat-induced stress on thermal pain. \cite{Attar2016-hi} investigated the thermo-visco-elastic behavior of tumorous and healthy bovine liver tissue treating the tissue as viscoelastic.  They have performed experimental tests and compared them with numerical simulations. 
Under specific power dissipation conditions, the results of the experimental tests are found to be in good agreement with results of the numerical solution, with experimental results slightly under performing the numerical solution.

We are proposing and investigating a hydrodynamic thermo-elasticity system in the context of mixture models where the mixture consists of a fluid and a solid phase.
This system is linear with the exception of the heat advection in the fluid phase.
Our main motivating application is biological tissues especially related to tumour growth; although at this time our model does not include any description of this growth.
As a first step, we aim to understand better the interplay of thermal coupling with the hydrodynamics on short-time scales where the tumour is essentially stable.
This follows the same line of reasoning as in
 \cite{alam2019mathematical,alam2022existence}.
We point out, that this type of model also comes up in porous media more generally so it is not restricted to this application.

We show that this coupled system has a unique weak solution, where the main difficulty is navigating the different coupling mechanisms in the system.
This analysis is done by combining standard existence results for elliptic PDEs with the theory for implicit parabolic PDEs (\cite{S96}).
Moreover, we analyze the system using a simplified 1D Cartesian coordinate framework and consider a small timescale.
On the short timescale considered here, the influence of the solid velocity on the fluid dynamics is neglected, which leads to a one-way coupling of the elastohydrodynamic equations.
However, the thermodynamic governing equations remain two-way coupled. Focusing on a one-dimensional setup reduces the complexity allowing for semi-analytical approach to analyze the underlying elastothermodynamics.
As a next step, we plan to further investigate this system by looking at larger time scales where both the growth of the tumour as well as some nonlinear material relationships have to be considered.

This manuscript is structured as follows: In \cref{sec:model}, we introduce the mathematical model we are considering, perform a non-dimensionalization, and introduce the corresponding weak forms.
This is followed by \cref{sec:analysis}, where we tend to the analysis of the weak form and show that there is a unique weak solution.
In \cref{sec:exact_solutions}, we calculate and present exact solutions in some simplified scenarios (1D Cartesian coordinates, isotropic permeability, time-independent) and go in \cref{sec:results_discussion} to discuss some rudimentary numerical results mainly aimed at highlighting how the different involved mechanism interact with each other.
%%%%%%%%%%%%%%%%%%%%%%%%%%%%%%%%%%%%%%%%%%%%%%%%%%%%%%%%%%%%%%%%%
% Problem formulation
%%%%%%%%%%%%%%%%%%%%%%%%%%%%%%%%%%%%%%%%%%%%%%%%%%%%%%%%%%%%%%%%%
\section{Mathematical Model and Setting}\label{sec:model}
In this section, we provide the mathematical model that we are considering in this work and present mathematical framework in which we analyze this model.
In the following, let $\Omega\subset\mathbb{R}^d$,  $d=1,2,3$, denote a bounded Lipschitz domain with outer unit normal vector $n$ that represents the tumor region and let $S=(0,T)$ be the time interval of interest.
We introduce a partition of the boundary $\Gamma,\Gamma_N\subset\partial\Omega$ such that $|\Gamma|>0$ and $\Gamma\cup\Gamma_N=\partial\Omega$.

We consider a bi-phasic mixture model for a porous tumor comprising of a solid phase (tumour tissue) and a liquid phase (extracellular fluid) where the main assumption is that both phases are present at every point.
For more details regarding this general modeling approach, we refer to \cite{B76,T84}, but in short: all relevant quantities, e.g., fluid and solid mass densities, are assumed to be defined in all of $\Omega$ and the exchange of momentum, mass, and heat between the two phases is accounted for via additional exchange terms in the corresponding balance equations.
This is an established starting point for the mathematical modelling of, in particular, tissue mechanics, see, e.g., \cite{BP03,DS18}.
In the following, we denote with subscripts $f$ and $s$ the various physical variables corresponding to the fluid phase and solid phase respectively.
For multi-phase systems, we refer to \cite{Ateshian2007,Giverso2015}.

% \begin{table}[h]
% \footnotesize
% \centering
% \begin{tabular}{l"l"l"}
% Symbol  						& Quantity 										& Units	 \\ \tline
% $\rho_i$						& Mass densities 							& \si{\kilogram\per\meter\cubed}     \\ 
% $\mu_i$							& 2nd Lamé parameter 					& \si{\kilogram\per\meter\per\second\squared}      \\
% $\lambda_f,\chi_s$	&	1st Lamé parameter					& \si{\kilogram\per\meter\per\second\squared}\\
% $\varphi_i$ 				& Volume fractions						& 1	\\
% $K$									&	Drag coefficient						& \si{\meter\squared\per\second}\\
% $b_i$								& Body force densities				&	\si{\kilogram\per\meter\squared\per\second\squared}\\
% $S_f(P)$								&	Fluid phase generation			&	\si{\per\second}\\[0.3cm]\hline
% $c_i$								& Specific heat capacity			& \si{\meter\squared\per\kelvin\per\second\squared}\\
% $\kappa_i$					& Heat conductivity						&	\si{\kilogram\meter\per\second\cubed\per\kelvin}\\
% $\alpha$						& Heat exchange coefficient 	&\si{\kilogram\per\meter\per\second\cubed\per\kelvin}\\
% $\gamma_i$					&	Modulus of heat dissipation	&	\si{\kilogram\per\meter\per\second\squared}\\
% $g_i$								& Heat production density			&	\si{\kilogram\per\meter\per\second\cubed}
% \end{tabular}
% \caption{List of model parameters with symbols and units}
% \end{table}

\begin{table}[h]
\footnotesize
\centering
\begin{tabular}{l"l"l"}
Symbol & Quantity & Units \\ \tline

$\rho_i$
& Partial mass density
& \si{\kilogram\per\meter\cubed}
\\

$\mu_f$
& Dynamic shear viscosity of the fluid
& \si{\kilogram\per\meter\per\second}
\\

$\lambda_f$
& Bulk viscosity coefficient of the fluid
& \si{\kilogram\per\meter\per\second}
\\

$\mu_s$
& Shear modulus of the solid
& \si{\kilogram\per\meter\per\second\squared}
\\

$\chi_s$
& First Lamé parameter of the solid
& \si{\kilogram\per\meter\per\second\squared}
\\

$\varphi_i$
& Volume fraction
& 1
\\

$K$
& Hydraulic mobility
& \si{\meter\cubed\second\per\kilogram}
\\

$b_i$
& Body force density
& \si{\kilogram\per\meter\squared\per\second\squared}
\\

$S_f(P)$
& Fluid volume source
& \si{\per\second}
\\

$a_0$
& Pressure-dependent source coefficient
& \si{\meter\second\per\kilogram}
\\

$P_a$
& Ambient vascular pressure
& \si{\kilogram\per\meter\per\second\squared}
\\[0.3cm]\hline

$c_i$
& Specific heat capacity
& \si{\meter\squared\per\kelvin\per\second\squared}
\\

$\kappa_i$
& Thermal conductivity
& \si{\kilogram\meter\per\second\cubed\per\kelvin}
\\

$h$
& Interphase heat exchange coefficient
& \si{\kilogram\per\meter\per\second\cubed\per\kelvin}
\\

$\alpha_s$
& Solid thermal stress coefficient
& \si{\kilogram\per\meter\per\second\squared\per\kelvin}
\\

$\gamma_i$
& Thermo-mechanical coupling coefficient
& \si{\kilogram\per\meter\per\second\squared}
\\

$g_i$
& Heat production density
& \si{\kilogram\per\meter\per\second\cubed}

\end{tabular}
\caption{List of model parameters with symbols and units.}
\end{table}

Now, let $\mathbf{v}_f$ and $\mathbf{v}_s$ denote the interstitial fluid velocity (IFV) and the velocity of the solid phase (SPV), respectively.
In addition, we introduce the volume fractions $\varphi_i$ $(i=s,f)$, the intrinsic densities $\widetilde\rho_i$, and the partial mass densities
$\rho_i=\phi_i\widetilde\rho_i$.
The mass balance equations are then given via
\[
\partial_t\rho_i+\nabla\cdot(\rho_i\mathbf{v}_i)=\tilde{\rho_i}S_i,
\]
where the $S_i$ account for the potential mass exchange between the two phases (as a consequence of, e.g., adsorption processes).
On the short time scale considered here, we do not resolve the individual evolution of the phase volume fractions.
Instead, we freeze them at constant reference values
\[
    \varphi_f,\varphi_s\in(0,1),
    \qquad
    \varphi_f+\varphi_s=1,
\]
and retain only the reduced mixture volume balance derived below.
Accordingly, the individual phase mass balances are used here only as a motivation for the reduced balance law and are not imposed as additional equations in the final model.
Similarly, the momentum balance equations can be stated as
\begin{align*}
\rho_i\left(\frac{\partial \mathbf{v}_i}{\partial t}+(\nabla \mathbf{v}_i)\mathbf{v}_i\right)&=\nabla\cdot \mathbf{T}_i+\mathbf{b}_i+{\bf \Pi}_i.
\end{align*}
Here, $\mathbf{T}_i$ are the Cauchy stress tensors, $\mathbf{b}_i$ volume force densities (modeling, e.g., gravity), and ${\bf \Pi}_i$ account for the exchange of momentum between the phases.
By Newton's third law this necessitates ${\bf \Pi}_f=-{\bf \Pi}_s$.
Moreover, postulating angular momentum balance in a non-polar setting (in particular, no exchange of angular momentum between the phases), we additionally get symmetry of stress tensors, i.e., ${\bf T}_i={\bf T}_i^T$. 
For the more complex situation allowing for angular momentum coupling, we refer to \cite{D85}.
Finally, the balance of internal energy for each phase is written in conservative form as
\[
    \partial_t(\rho_i e_i)
    + \nabla\cdot(\rho_i e_i \mathbf v_i + \mathbf q_i)
    = h_i + \Sigma_i,
\]
where $e_i$ denotes the specific internal energy, $\mathbf q_i$ the conductive
heat flux, $h_i$ the volumetric heat source, and $\Sigma_i$ accounts for heat
exchange and thermo-mechanical coupling.

These abstract balance equations, which are in principle valid in a wide range of situations, have to be completed with constitutive relations adapted to the particular setting we have in mind.
We start by introducing some additional variables, namely the solid deformation $\mathbf{u}_s$ which satisfies $\partial_t\mathbf{u}_s=\mathbf{v}_s$,  the pressure $P$, and temperatures $\theta_f$ and $\theta_s$.
We work with the overall pressure $P$ (instead of separate phase pressures) as is usually done in mixture theory, e.g., \cite{AP02, P96}.
The system we consider in this work is formulated with respect to the state variables $(\mathbf{v}_f, \mathbf{u}_s, P, \theta_f, \theta_s)$.
In some cases, the evolution of the volume fractions $\varphi_i$ is also taken into account as part of the overall problem (e.g., \cite{BP03,DS16}) but this is not the case here.

Assuming the intrinsic phase densities $\widetilde\rho_i$ to be constant in space and time,
the mass balance equations can equivalently be written in terms of the volume fractions as
\[
    \partial_t\phi_i+\nabla\cdot(\phi_i v_i)=S_i.
\]
We furthermore assume the mixture to remain saturated and use constant
reference volume fractions satisfying
\[
    \varphi_f+\varphi_s=1.
\]
Formally adding the phase balances and then applying the frozen-volume-
fraction approximation gives
\[
    \nabla\cdot
    \left(
        \varphi_f\mathbf v_f
        +
        \varphi_s\partial_t\mathbf u_s
    \right)
    =
    S_f+S_s.
\]
In a closed system one would have $S_f+S_s=0$, but this assumption is not appropriate in the present setting.
Since the time scale of tumour growth is significantly larger than that of the interstitial fluid processes considered here
\cite{BJ89, NB97}, we set $S_s=0$.
The fluid source is assumed to be driven by the average transmural pressure and is modeled by
\[
    S_f=a_0(P_a-P),
\]
where $P_a$ denotes the ambient vascular pressure; see, e.g., \cite{DS16,NB97}.
On the short time scale considered here, the velocity of the solid phase is assumed to be small compared with the interstitial fluid velocity.
We therefore neglect the solid-velocity contribution in the mixture volume balance and obtain
\[
    \nabla\cdot(\phi_f v_f)=a_0(P_a-P).
\]
The same approximation is used in the interphase momentum exchange below.

We assume that the fluid phase behaves like a \emph{Newtonian} fluid and that the solid stresses can be described with linear thermoelasticity coupled with an additional pressure stress due to the fluid phase:
\begin{align}
\mathbf{T}_f&=-\varphi_fP\mathds{I}+\lambda_f\nabla\cdot \mathbf{v}_f\mathds{I}+2\mu_fD(\mathbf{v}_f),\\
\mathbf{T}_s&=-\varphi_sP\mathds{I}+\chi_s\nabla\cdot \mathbf{u}_s\mathds{I}+2\mu_sD(\mathbf{u}_s)-\alpha_s\theta_s\mathds{I}.
\end{align}
Here, $D(\mathbf{v})=\nicefrac{1}{2}(\nabla \mathbf{v}+\nabla \mathbf{v}^T)$ denotes the symmetrized gradient, $\lambda_f$ and $\mu_f$ the first and second coefficients of viscosity, and $\chi$ and $\mu_s$ the elasticity coefficients.
Finally, $\alpha_s$ is the thermal stress coefficient of the solid phase.
Similar constitutive laws without the thermal part can be found in \cite{AP02, BP03}; for the linear thermoelasticity, we refer to \cite{K79}.
The momentum exchange is assumed to be driven primarily by the difference in velocities (see, e.g., \cite{BP03,DS16,AP02} for more details), that is
\[
{\bf \Pi}_f=-{\bf \Pi}_s=\frac{1}{K}(\partial_t\mathbf{u}_s-\mathbf{v}_f)
\]
with $\nicefrac{1}{K}$ denoting the drag coefficient.
The solid is generally assumed to be close to rigid and its velocity therefore negligible compared to the fluid leading to ${\bf \Pi}_f=-{\bf \Pi}_s=-K^{-1}\mathbf{v}_f$ \cite{AP02,Rajagopal2007,Turner2015}; note that this is consistent with the short-time approximation introduced above.
In addition, inertial effects are neglected in both phases, so that the momentum balances are treated quasi-statically.
Consequently, the momentum balance equations for the two phases take the form of a coupled system of quasi-stationary elliptic PDEs:
\begin{subequations}\label{full_model}
\begin{alignat}{2}
-\nabla\cdot\left(2\mu_fD(\mathbf{v}_f)+\lambda_f\nabla\cdot \mathbf{v}_f\mathds{I}-\varphi_fP\mathds{I}\right)+\frac{1}{K}\mathbf{v}_f&=\mathbf{b}_f&\quad&\text{in}\ S\times\Omega,\label{full_model:1}\\
-\nabla\cdot\left(2\mu_sD(\mathbf{u}_s)+\chi_s\nabla\cdot \mathbf{u}_s\mathds{I}-\varphi_sP\mathds{I}-\alpha_s\theta_s\mathds{I}\right)-\frac{1}{K}\mathbf{v}_f&=\mathbf{b}_s&\quad&\text{in}\ S\times\Omega,\label{full_model:1b}\\
\nabla\cdot(\varphi_f\mathbf{v}_f)+a_0P&=a_0P_a&\quad&\text{in}\ S\times\Omega.\label{full_model:2}
\end{alignat}
% \hl{Note: We need to consider the following mass balance equation if} \(\mathbf{v}_s\neq 0\).
% \[{\nabla\cdot(\varphi_f\mathbf{v}_f+\varphi_s\mathbf{v}_s)=a_0(P_a-P)}\quad\text{in}\ S\times\Omega.\]\\
%
At the outer boundary, we assume the solid part to be clamped and the fluid phase to obey a non-slip condition at one part of the boundary mixed with a (potentially non-homogeneous) Neumann boundary condition on the remaining part:
\begin{align}
    \mathbf{u}_s&=0\quad\text{on}\ S\times\partial\Omega,\label{full_model:3}\\
    \mathbf{v}_f&=0\quad\text{on}\ S\times\Gamma,\qquad\label{full_model:4}\\
    (2\mu_fD(\mathbf{v}_f)+\lambda_f\nabla\cdot \mathbf{v}_f\mathds{I}-\varphi_fP\mathds{I})\cdot{\bf n}&={\bf g}\quad\text{on}\ S\times\Gamma_N.\label{full_model:5}
\end{align}
Here, $g\colon S\times\Gamma_N\to\R^3$ denotes potential forces acting on the fluid via the surface $\Gamma_N$.

The heat model considered in this work is based on the principles of linear thermoelasticity, structurally identical to poroelasticity \cite{K79,SM02}.
Here, the internal energies are linear functions of the temperatures, i.e.,
\[
e_i=c_i\theta_i
\]
where $c_i>0$ denote the specific heat capacities of the phases.
The conductive heat fluxes are given by Fourier's law,
\[
q_i=-\kappa_i\nabla\theta_i
\]
with conductivity matrices $\kappa_i\in\R_{\rm sym}^{d\times d}$.
Here, we neglect convection in the solid phase consistent with the small-solid-velocity approximation introduced above.
Assuming a linear heat exchange law of Newton type, i.e., $h(\theta_f-\theta_s)$, as well as a linear thermo-mechanical coupling $(\gamma_i\nabla\cdot \mathbf{v}_i)$, we end up with
\begin{alignat}{2}
\rho_fc_f\partial_t\theta_f-\kappa_f\nabla^2\theta_f+\rho_fc_f\nabla\cdot(\theta_f\mathbf{v}_f)+h(\theta_f-\theta_s)&=
-\gamma_f\nabla\cdot\mathbf{v}_f+g_f&\quad&\text{in}\ S\times\Omega,\label{full_model:6}\\
\rho_sc_s\partial_t\theta_s-\kappa_s\nabla^2\theta_s-h(\theta_f-\theta_s)&=
-\gamma_s\nabla\cdot\partial_t\mathbf{u}_s+g_s&\quad&\text{in}\ S\times\Omega.\label{full_model:7}
\end{alignat}
The parameter $h$ is the heat exchange coefficient and $\gamma_i$ are the thermo-mechanical coupling coefficients.
Finally, $g_i$ denotes the heat production densities.
We complete the system with the initial conditions
\begin{equation}
\theta_f(0)=\theta_s(0)=\theta_0\quad\text{in}\ \Omega,\label{full_model:8}
\end{equation}
and Dirichlet boundary conditions
\begin{equation}
\theta_f=\theta_s=\theta_a\quad \text{on}\ S\times\partial\Omega\label{full_model:9}
\end{equation}
\end{subequations}
where $\theta_a$ denotes the ambient temperature.
Of course, other types of boundary conditions, in particular Neumann and Robin exchange conditions, are also reasonable here.
The mathematical model describing the hydrodynamics with thermo-elastic effects is then given by \cref{full_model:1,full_model:1b,full_model:2,full_model:3,full_model:4,full_model:5,full_model:6,full_model:7,full_model:8,full_model:9}.

\subsection{Nondimensionalization}
We introduce the dimensionless variables
\[
    x=L\hat x,\qquad
    t=\frac{L}{V}\hat t,\qquad
    \mathbf v_f=V\hat{\mathbf v}_f,\qquad
    \mathbf u_s=\frac{\mu_fV}{\mu_s}\hat{\mathbf u}_s,
\]
\[
    P=\frac{\mu_fV}{L}\hat P,\qquad
    \theta_i-\theta_a=W_i\hat\theta_i,\qquad
    \mathbf b_i=\frac{\mu_fV}{L^2}\hat{\mathbf b}_i,\qquad
    g_i=\frac{\kappa_iW_i}{L^2}\hat g_i,
    \qquad i=f,s.
\]
Moreover, the boundary traction is scaled according to
\[
    \mathbf g=\frac{\mu_fV}{L}\hat{\mathbf g}.
\]
Due to the shift by the ambient temperature in $\theta_i-\theta_a=W_i\hat\theta_i$, the convective term in the fluid heat equation transforms according to
\[
    \rho_f c_f \nabla\cdot(\theta_f\mathbf v_f)
    =
    \frac{\rho_f c_f V W_f}{L}
    \nabla\cdot(\widehat\theta_f\widehat{\mathbf v}_f)
    +
    \frac{\rho_f c_f V\theta_a}{L}
    \nabla\cdot\widehat{\mathbf v}_f.
\]
We therefore combine the latter contribution with the
thermo-mechanical coupling term and define
\[
    \delta'_f
    :=
    \frac{LV}{\kappa_fW_f}
    \left(
        \gamma_f+\rho_fc_f\theta_a
    \right).
\]
Here, $L$, $V$, and $W_i$ denote characteristic length, fluid velocity,
and phase-specific temperature scales, respectively. 
For notational convenience, we drop the hats in the following and denote the non-dimensional body forces and boundary traction by
\[
    f_f:=b_f,
    \qquad
    f_s:=b_s,
    \qquad
    f_N:=g.
\]
The non-dimensional hydromechanical system is given by
\begin{subequations}
\begin{alignat}{2}
-\nabla\cdot\Big(
2D(\mathbf v_f)
+\big(\Lambda_f\nabla\cdot\mathbf v_f-\varphi_f P\big)\mathds I
\Big)
+\frac{1}{Da}\mathbf v_f
&=f_f
&&\quad\text{in }S\times\Omega,
\label{main-5a}
\\
-\nabla\cdot\Big(
2D(\mathbf u_s)
+\big(\Lambda_s\nabla\cdot\mathbf u_s-\varphi_s P-\delta_s\theta_s\big)\mathds I
\Big)
-\frac{1}{Da}\mathbf v_f
&=f_s
&&\quad\text{in }S\times\Omega,\label{main-5b}
\\
\nabla\cdot(\varphi_f\mathbf v_f)+a_1P
&=a_2
&&\quad\text{in }S\times\Omega,\label{main-5c}
\end{alignat}
with boundary conditions
\begin{align}
\mathbf u_s&=0
&&\quad\text{on }S\times\partial\Omega,\label{main-5d}
\\
\mathbf v_f&=0
&&\quad\text{on }S\times\Gamma,\label{main-5e}
\\
\Big(
2D(\mathbf v_f)
+\big(\Lambda_f\nabla\cdot\mathbf v_f-\varphi_fP\big)\mathds I
\Big)\mathbf n
&=f_N
&&\quad\text{on }S\times\Gamma_N.\label{main-5f}
\end{align}
The corresponding non-dimensional thermal equations read
\begin{alignat}{2}
\mathrm{Pe}_f
\Big(
\partial_t\theta_f
+\nabla\cdot(\theta_f\mathbf v_f)
\Big)
-\nabla^2\theta_f
&=
-N(\theta_f-W\theta_s)
-\delta'_f\nabla\cdot\mathbf v_f
+g_f
&&\quad\text{in }S\times\Omega,
\label{main-5g}
\\
\mathrm{Pe}_s\partial_t\theta_s
-\nabla^2\theta_s
&=
-\kappa N
\left(
\theta_s-\frac{1}{W}\theta_f
\right)
-\delta'_s\nabla\cdot\partial_t\mathbf u_s
+g_s
&&\quad\text{in }S\times\Omega.
\label{main-5h}
\end{alignat}
Finally, the initial and boundary conditions are
\begin{align}
\theta_f(0)&=\theta_{f,0},
&
\theta_s(0)&=\theta_{s,0}
&&\quad\text{in }\Omega,\label{main-5i}
\\
\theta_f&=0,&\theta_s&=0
&&\quad\text{on }S\times\partial\Omega.\label{main-5j}
\end{align}
\end{subequations}

\begin{table}[h]
\footnotesize
\centering
\begin{tabular}{|c|c|}
\hline
\textbf{Non-dimensional parameters} & \textbf{Description} \\
\hline

$Da=\dfrac{\mu_f K}{L^2}$
& Darcy number \\
\hline

$a_1=a_0\mu_f,\qquad
 a_2=\dfrac{a_0P_aL}{V}$
& Non-dimensional pressure-source parameters \\
\hline

$\Lambda_f=\dfrac{\lambda_f}{\mu_f},\qquad
 \Lambda_s=\dfrac{\chi_s}{\mu_s}$
& Fluid viscosity ratio and solid elasticity ratio \\
\hline

$\delta_s=\dfrac{L\alpha_sW_s}{V\mu_f}$
& Ratio of thermoelastic stress to viscous stress \\
\hline

$\mathrm{Pe}_f=\dfrac{\rho_fc_fVL}{\kappa_f}$
& Péclet number of the fluid phase \\
\hline

$\mathrm{Pe}_s=\dfrac{\rho_sc_sVL}{\kappa_s}$
& Péclet number of the solid phase \\
\hline

$N=\dfrac{hL^2}{\kappa_f}$
& Non-dimensional interphase heat exchange parameter \\
\hline

$\kappa=\dfrac{\kappa_f}{\kappa_s},\qquad
 W=\dfrac{W_s}{W_f}$
& Ratios of thermal conductivities and temperature scales \\
\hline

$\displaystyle
\delta'_f=
\frac{LV(\gamma_f+\rho_fc_f\theta_a)}
     {\kappa_fW_f}
$
& Fluid effective divergence-coupling parameter \\
\hline
$
\delta'_s=
\frac{\gamma_s\mu_fV^2}
     {\mu_s\kappa_sW_s}
$
& Solid
thermo-mechanical coupling parameter \\
\hline

\end{tabular}
\caption{\label{non-dim}List of non-dimensional parameters.}
\end{table}

\subsection{Variational formulation}
We introduce the Hilbert spaces 
\begin{align*}
\mathbf{H}^1_{0}(\Omega)&:=H_0^1(\Omega)^d,\\
\mathbf{H}^1_{\Gamma}(\Omega)&:=\{u\in H^1(\Omega)^d \ : u=0\ \text{on}\ \Gamma\},\\
\mathcal{W}(S,\Omega)&:=\{\theta\in L^2(S;H^1_0(\Omega)) \ : \ \partial_t\theta\in L^2(S;H^{-1}(\Omega))\}.
\end{align*}
Please note that Korn's inequality \cite[Theorem 3.1]{Duvaut1976} together with Poincare's inequality ensures that there is $C>0$ such that
\begin{equation}\label{korn}
   \|u\|_{\mathbf{H}^1_{0}(\Omega)}\leq C\|D(u)\|_{L^2(\Omega)},\quad \|v\|_{\mathbf{H}^1_{\Gamma}(\Omega)}\leq C\|D(v)\|_{L^2(\Omega)}.
\end{equation}
In the following, we use $\langle\cdot,\cdot\rangle_{V}$ to denote the dual product between $V$ and its dual $V^*$.
A weak form of the hydromechanics part in its nondimensional form, \cref{main-5a,main-5b,main-5c,main-5d,main-5e,main-5f}, is given via
\begin{subequations}
\begin{align}
\int_\Omega \mathcal{A}_fD(\mathbf{v}_f):D(\mathbf{w}_f)-\varphi_fP\nabla\cdot \mathbf{w}_f+\frac{1}{Da}\mathbf{v}_f\cdot \mathbf{w}_f\di{x}=\int_\Omega\hspace{-.1cm} f_f\cdot \mathbf{w}_f\di{x}
+\int_{\Gamma_N}\hspace{-.3cm}f_N\cdot \mathbf{w}_f\di{\sigma}\label{weak1},\\
\int_\Omega \mathcal{A}_sD(\mathbf{u}_s):D(\mathbf{w}_s)-\left(\varphi_sP+\delta_s\theta_s\right)\nabla\cdot \mathbf{w}_s-\frac{1}{Da}\mathbf{v}_f\cdot \mathbf{w}_s\di{x}=\int_\Omega f_s\cdot \mathbf{w}_s\di{x}\label{weak2},\\
\int_\Omega\left(\varphi_f\nabla\cdot \mathbf{v}_f-a_2+a_1P\right)Q\di{x}=0,\label{weak3}
\end{align}
which is assumed to hold for almost all $t\in S$ and for all test functions
\[
(\mathbf{w}_f,\mathbf{w}_s,Q)
\in
H^1_\Gamma(\Omega)^d
\times
H^1_0(\Omega)^d
\times
L^2(\Omega).
\]
For ease of notation, we have introduced the fourth-rank tensors $ A_i\in\R^{d\times d\times d\times d}$ via\footnote{Here, $\delta_{ij}$ is the Kromecker delta.}
\[
( A_i)_{klmn}
=
\delta_{km}\delta_{ln}
+\delta_{kn}\delta_{lm}
+\Lambda_i\delta_{kl}\delta_{mn},
\]
The variational form of the thermal subsystem in its non-dimensionalized form \cref{main-5g,main-5h,main-5i,main-5j} is given by
\begin{multline}
-\int_S\int_\Omega
\mathrm{Pe}_f\,\theta_f\,\partial_t\psi_f
\di{x}\di{t}
+\int_S\int_\Omega
\delta'_f\,\nabla\cdot\mathbf{v}_f\,\psi_f
\di{x}\di{t}
\\
+\int_S\int_\Omega
\nabla\theta_f\cdot\nabla\psi_f
\di{x}\di{t}
-\int_S\int_\Omega
\mathrm{Pe}_f\,\theta_f\mathbf{v}_f\cdot\nabla\psi_f
\di{x}\di{t}
\\
+\int_S\int_\Omega
N(\theta_f-W\theta_s)\psi_f
\di{x}\di{t}
=
\int_S\int_\Omega
g_f\psi_f
\di{x}\di{t}
+\int_\Omega
\mathrm{Pe}_f\,\theta_{f,0}\psi_f(0)
\di{x},
\label{weak4}
\end{multline}
and
\begin{multline}
-\int_S\int_\Omega
\mathrm{Pe}_s\,\theta_s\,\partial_t\psi_s
\di{x}\di{t}
-\int_S\int_\Omega
\delta'_s\,\nabla\cdot\mathbf{u}_s\,\partial_t\psi_s
\di{x}\di{t}
\\
+\int_S\int_\Omega
\nabla\theta_s\cdot\nabla\psi_s
\di{x}\di{t}
+\int_S\int_\Omega
\kappa N
\left(
\theta_s-\frac{1}{W}\theta_f
\right)\psi_s
\di{x}\di{t}
\\
=
\int_S\int_\Omega
g_s\psi_s
\di{x}\di{t}
+\int_\Omega
\left(
\mathrm{Pe}_s\theta_{s,0}
+\delta'_s\nabla\cdot\mathbf{u}_s(0)
\right)\psi_s(0)
\di{x}.
\label{weak5}
\end{multline}
These identities are required to hold for all test functions
\[
    \psi_f,\psi_s\in \left\{
        \psi\in L^2(S;H_0^1(\Omega))
        \cap H^1(S;L^2(\Omega))
        :
        \psi(T)=0
    \right\}.
\]
\end{subequations}
One potential issue in \cref{weak5} concerns the time regularity of
$\mathbf u_s$.
Since \cref{weak2} is only required to hold for almost
all $t\in S$, it is a priori unclear how to interpret
$\nabla\cdot\mathbf u_s(0)$.
To ensure that this initial value is well-defined, we impose additional
time regularity on the hydrodynamic data $f_s$, $f_f$, and $f_N$.

\begin{remark}
A less restrictive notion of weak solution could alternatively be used,
in which the initial condition is imposed on the combined quantity
\[
    \lim_{t\to0}
    \left(
        \mathrm{Pe}_s\theta_s(t)
        +\delta'_s\nabla\cdot\mathbf u_s(t)
    \right).
\]
In this setting, neither $\theta_s$ nor
$\nabla\cdot\mathbf u_s$ needs to possess a weak time derivative
separately; it is sufficient to require
\[
    \partial_t
    \left(
        \mathrm{Pe}_s\theta_s
        +\delta'_s\nabla\cdot\mathbf u_s
    \right)
    \in L^2(S;H^{-1}(\Omega)).
\]
We refer to \cite{SM02}, where this approach is explored in more detail
in the context of poroelasticity.
\end{remark}

We provide a short summary of the main regularity assumptions:
\begin{itemize}
    \item[$\mathbf{(A1)}$]\label{ass1}\textbf{Geometry:} $\Omega\subset\R^d$ is a bounded Lipschitz domain with outer normal vector $n$.
    Its boundary is decomposed in two disjunct components $\Gamma,\Gamma_N$ with $\Gamma_N\cup\Gamma=\partial\Omega$ and positive surface measures $|\Gamma_N|,|\Gamma|>0$.
    \item[$\mathbf{(A2)}$]\label{ass2}\textbf{Data:} $f_f\in L^\infty(S;L^2(\Omega)^d)$,
    $f_s\in L^2(S\times\Omega)^d$,
    $g_f,g_s\in L^2(S\times\Omega)$,
    $f_N\in L^\infty(S;L^2(\Gamma_N)^d)$,
    and
    $\theta_{f,0},\theta_{s,0}\in L^2(\Omega)$.
    \item[$\mathbf{(A3)}$]\label{ass3}\textbf{Time-regularity of data:} $\partial_tf_f\in L^2(S;(\mathbf{H}_{\Gamma}^{1}(\Omega))^*)$, $\partial_tf_N\in L^2(S\times\Gamma_N)^d$, and $\partial_tf_s\in L^2(S;\mathbf{H}^{-1}(\Omega))$
    \item[$\mathbf{(A4)}$]\label{ass4}\textbf{Coefficients:} The parameters
    \[
        Da,\ a_1,\ \Lambda_f,\ \Lambda_s,\ 
        \mathrm{Pe}_f,\ \mathrm{Pe}_s,\ N,\ \kappa,\ W
    \]
    are strictly positive, while
    \[
        \varphi_f,\varphi_s\in(0,1),
        \qquad
        \varphi_f+\varphi_s=1.
    \]
    Moreover, $\delta_s,\delta'_f,\delta'_s\geq0$.
\end{itemize}

%%%%%%%%%%%%%%%%%%%%%%%%%%%%%%%%%%%%%%%%%%%%%%%%%%%%%%%%%%%%%%%%%%%%
%%%%%%%%%%%%%%%%%%%%%%%%%%%%%%%%%%%%%%%%%%%%%%%%%%%%%%%%%%%%%%%%%%%%
%%%%%%%% Analysis
%%%%%%%%%%%%%%%%%%%%%%%%%%%%%%%%%%%%%%%%%%%%%%%%%%%%%%%%%%%%%%%%%%%%
%%%%%%%%%%%%%%%%%%%%%%%%%%%%%%%%%%%%%%%%%%%%%%%%%%%%%%%%%%%%%%%%%%%%

\section{Well-posedness}\label{sec:analysis}
In this section, we tend to the mathematical analysis of the coupled problem with quasi-static hydrodynamics in its non-dimensional form as given by \cref{weak1,weak2,weak3,weak4,weak5}.

Our approach here is partially motivated by \cite{SM02} where a somewhat similar structure to our thermo-elasticity coupling is considered in the context of poroelasticity, although without the mixture component and with different coupling mechanisms.
We also refer to \cite{Eden_Muntean2017} where the analysis of a thermo-elasticity model is considered.
Regarding the general strategy we employ, the analysis of this coupled system is done in two parts:
\begin{itemize}
    \item[1)] We first look at the hydrodynamics sub-problem given by \cref{weak1,weak2,weak3} for prescribed temperature profiles.
    This is done in the context of standard elliptic theory.
    The main result here is given by \cref{lemma_step1}.
    \item[2)] Using the resulting solution operator, we tackle the complete coupled problem.
    In this part, we are relying on results on implicit partial differential equations as presented in \cite{S96}.
    See \cref{theorem:existence}.
\end{itemize}

We first focus on the hydrodynamics part given by \cref{weak1,weak2,weak3} with some prescribed temperature function $\vartheta_s\in L^2(\Omega)$.
Making explicit use of \cref{weak3}, we see that 
\begin{equation}\label{eq:pressure}
    P=\frac{1}{a_1}
    \left(
        a_2-\varphi_f\nabla\cdot\mathbf v_f
    \right)
    \qquad\text{in }L^2(\Omega).
\end{equation}
Substituting \cref{eq:pressure} into \cref{weak1,weak2}, we are therefore able to decouple the pressure function from the fluid velocity.
By noting that 
\[
\int_\Omega\nabla\cdot \mathbf{w}_f\di{x}=\int_{\Gamma_N}\mathbf{w}_f\cdot n\di{\sigma} \quad (\mathbf{w}_f\in\mathbf{H}^1_{\Gamma}(\Omega)),\qquad\int_\Omega\nabla\cdot \mathbf{w}_s\di{x}=0 \quad (\mathbf{w}_s\in\mathbf{H}^1_{0}(\Omega)),
\]
we arrive at the following problem:
Find $(\mathbf{v}_f,\mathbf{u}_s)\in \mathbf{H}^1_{\Gamma}(\Omega)\times \mathbf{H}^1_{0}(\Omega)$ such that
\begin{subequations}
\begin{multline}
\int_\Omega
\mathcal A_f D(\mathbf v_f):D(\mathbf w_f)
+
\frac{\varphi_f^2}{a_1}
\nabla\cdot\mathbf v_f\,
\nabla\cdot\mathbf w_f
+
\frac{1}{Da}\mathbf v_f\cdot\mathbf w_f
\di{x}\\
=
\int_\Omega f_f\cdot\mathbf w_f\di{x}
+
\int_{\Gamma_N}
\left(
    f_N
    +
    \frac{\varphi_fa_2}{a_1}\mathbf n
\right)
\cdot\mathbf w_f\di{\sigma},
\label{a:w4}
\end{multline}
\begin{align}
\int_\Omega
\mathcal A_s D(\mathbf u_s):D(\mathbf w_s)
+
\left(
    \frac{\varphi_s\varphi_f}{a_1}
    \nabla\cdot\mathbf v_f
    -
    \delta_s\theta_s
\right)
\nabla\cdot\mathbf w_s
-
\frac{1}{Da}\mathbf v_f\cdot\mathbf w_s
\di{x}
=
\int_\Omega f_s\cdot\mathbf w_s\di{x}.
\label{a:w5}
\end{align}
for all $(\mathbf{w}_f,\mathbf{w}_s)\in \mathbf{H}^1_{\Gamma}(\Omega)\times \mathbf{H}^1_{0}(\Omega)$.
\end{subequations}
From here, the pressure function can then be recovered via \cref{eq:pressure}.
For a more compact operator formulation of this system, we introduce elliptic operators
\begin{subequations}
\[
\mathcal{A}_f\colon \mathbf{H}^1_{\Gamma}(\Omega)\to (\mathbf{H}^1_{\Gamma}(\Omega))^*,\quad \mathcal{A}_s\colon \mathbf{H}^1_{0}(\Omega)\to \mathbf{H}^{-1}(\Omega)
\]
via 
\begin{align}
\langle\mathcal{A}_f\mathbf{v}_f, \mathbf{w}_f\rangle_{\mathbf{H}^1_{\Gamma}(\Omega)}&=\int_\Omega A_fD(\mathbf{v}_f):D(\mathbf{w}_f)\di{x}+\frac{\varphi_f^2}{a_1}\nabla\cdot \mathbf{v}_f\nabla\cdot \mathbf{w}_f\di{x},\label{operator1}\\
\langle \mathcal{A}_s\mathbf{u}_s,\mathbf{w}_s\rangle_{\mathbf{H}^1_{0}(\Omega)}&=\int_\Omega A_sD(\mathbf{u}_s):D(\mathbf{w}_s)\di{x}.\label{operator2}
\end{align}
Please note, that by Korn's and Poincaré's inequalities, the operators $A_f$ and $A_s$ are linear, continuous, and coercive on $H^1_\Gamma(\Omega)^d$ and $H_0^1(\Omega)^d$, respectively.
In addition, we introduce the coupling operator 
\[
\mathcal{B}\colon \mathbf{H}^1_{\Gamma}(\Omega)\to \mathbf{H}^{-1}(\Omega),\quad\nabla_s\colon L^2(\Omega)\to \mathbf{H}^{-1}(\Omega)
\]
defined via
\begin{align}
\langle \mathcal{B}\mathbf{v}_f,\mathbf{w}_s\rangle_{\mathbf{H}^1_{0}(\Omega)}&=\int_\Omega \frac{\varphi_s\varphi_f}{a_1}\nabla\cdot \mathbf{v}_f\nabla\cdot \mathbf{w}_s\di{x},\label{operator3}\\
\langle \nabla_s\pi,\mathbf{w}_s\rangle_{\mathbf{H}^1_{0}(\Omega)}&=-\int_\Omega \delta_s\pi\nabla\cdot \mathbf{w}_s\di{x}\label{operator4}
\end{align}
as well as the right-hand side operator
\begin{equation}\label{operator5}
\mathcal{F}_f\in (\mathbf{H}^1_{\Gamma}(\Omega))^*,\qquad \langle\mathcal{F}_f,\mathbf{w}_f\rangle_{\mathbf{H}^1_{\Gamma}(\Omega)}=\int_\Omega f_f\cdot \mathbf{w}_f\di{x}+\int_{\Gamma_N}(f_N+
    \frac{\varphi_fa_2}{a_1}\mathbf n)\cdot \mathbf{w}_f\di{\sigma}.
\end{equation}
\end{subequations}
With this, the weak form given by \cref{a:w4,a:w5} can alternatively be expressed as an abstract operator problem:
\begin{alignat*}{2}
\mathcal{A}_f\mathbf{v}_f+\frac{1}{Da}\mathbf v_f&=\mathcal{F}_f&\quad&\text{in}\ (\mathbf{H}^1_{\Gamma}(\Omega))^*,\\
\mathcal{A}_s\mathbf{u}_s+\left(\mathcal{B}-\frac{1}{Da}\ID\right)\mathbf v_f+\nabla_s\theta_s&=f_s&\quad&\text{in}\ \mathbf{H}^{-1}(\Omega).
\end{alignat*}
Here, $\ID$ denotes the identity operator.
In the following lemma, we establish that this system of linear elliptic PDEs has a unique solution.
%%%%%%%%%%%%%%%%%%%%%%%%%%%%%%%%%%%%%%%%%%%%%%%%
% Lemma: Hydro with fixed temperature
%%%%%%%%%%%%%%%%%%%%%%%%%%%%%%%%%%%%%%%%%%%%%%%%
\begin{lemma}\label{lemma_step1}
Let $f_f,f_s\in L^2(\Omega)$ and $f_N\in L^2(\Gamma_N)$.
For every $\theta_s\in L^2(\Omega)$, there is a unique solution $(\mathbf{v}_f,\mathbf{u}_s)\in \mathbf{H}^1_{\Gamma}(\Omega)\times \mathbf{H}^1_{0}(\Omega)$ to the problem given by \cref{a:w4,a:w5}.
The corresponding pressure function $P\in L^2(\Omega)$ is then given by \cref{eq:pressure}.
Moreover, the solution satisfies energy estimates of the form
\begin{subequations}\label{energy_estimates}
\begin{align}
\|\mathbf{v}_f\|_{\mathbf{H}^1_{\Gamma}(\Omega)}&\lesssim\|f_f\|_{L^2(\Omega)}+\|f_N\|_{L^2(\Gamma_N)}+1,\label{energy_estimates_a}\\
\|\mathbf{u}_s\|_{\mathbf{H}^1_{0}(\Omega)}&\lesssim\|\theta_s\|_{L^2(\Omega)}+\|f_f\|_{L^2(\Omega)}+\|f_s\|_{L^2(\Omega)}+\|f_N\|_{L^2(\Gamma_N)}+1\label{energy_estimates_b}.
\end{align}
\end{subequations}
\end{lemma}
\begin{proof}
For fixed $t\in S$, the fluid dynamics can be solved via the Lemma of Lax-Milgram (note that $(A_f+\frac1{Da}I)$ is linear, continuous, and coercive) such that
$$
\left(\mathcal{A}_f+\frac{1}{Da}I\right)\mathbf{v}_f=\mathcal{F}_f\quad\text{in}\ (\mathbf{H}^1_{\Gamma}(\Omega))^*,
$$
or, more explicitly,
\begin{equation}\label{eq:vf}
\mathbf{v}_f=\mathcal{S}_f\mathcal{F}_f\quad\text{in}\ \mathbf{H}^1_{\Gamma}(\Omega),
\end{equation}
where we have introduced the resolvent operator
\[
\mathcal{S}_f\colon\left(\mathbf{H}_\Gamma^1(\Omega)\right)^*\to \mathbf{H}_\Gamma^1(\Omega),\quad \mathcal{S}_f=\left(\mathcal{A}_f+\frac{1}{Da}I\right)^{-1}.
\]
With Korn's inequality, we have
\[
\langle\mathcal{A}_fW+\frac{1}{Da}\mathbf{W},\mathbf{W}\rangle_{\mathbf{H}_\Gamma^1(\Omega)}\ge c\|\nabla \mathbf{W}\|^2_{L^2(\Omega)}+\frac{1}{Da}\|\mathbf{W}\|^2_{L^2(\Omega)}
\]
for all $\mathbf{W}\in H^1_\Gamma(\Omega)$
which in turn implies
\[
\|\mathcal{S}_f\mathbf{W}\|_{L^2(\Omega)}\le Da\|\mathbf{W}\|_{L^2(\Omega)}
\]
for all $\mathbf W\in L^2(\Omega)$.
%
% We plug \cref{eq:vf} into the solid momentum equation and get
% \[
% \mathcal{A}_s\mathbf{u}_s+\nabla\theta_s=f_s-\left(\mathcal{B}-\frac{1}{Da}\right)\mathcal{R}\mathcal{F}_f.
% \]

For the solid phase, we similarly have linearity and continuity of $\mathcal{A}_s$.
Since it is also coercive by virtue of Korn's inequality \cref{korn} and since $\mathcal{B}$ is a linear and continuous operator, we find a unique $\mathbf{u}_s$ satisfying
\[
\mathcal{A}_s\mathbf{u}_s+\left(\mathcal{B}-\frac{1}{Da}\right)\mathbf{v}_f=-\nabla_s\theta_s+f_s\quad\text{in}\ \mathbf{H}^{-1}(\Omega),
\]
where $\mathbf{v}_f$ is given by \cref{eq:vf}.
Solving for $\mathbf{u}_s$ and using \cref{eq:vf}, we calculate
\begin{equation}\label{U_operator}
\mathbf{u}_s=-\mathcal{A}_s^{-1}\nabla_s\theta_s-\mathcal{A}_s^{-1}\left(\mathcal{B}-\frac1{Da}I\right)\mathcal{S}_f\mathcal F_f+\mathcal{A}_s^{-1}f_s\quad\text{in}\ \mathbf{H}^1_0(\Omega).
\end{equation}
We then introduce the heat-induced and data-induced displacements via
\[
u_\theta=-\mathcal A_s^{-1}\nabla_s\theta,\quad u_d=-\mathcal{A}_s^{-1}\left(\mathcal{B}-\frac1{Da}I\right)\mathcal{S}_f\mathcal F_f+\mathcal{A}_s^{-1}f_s
\]
and get
\[
\mathbf{u}_s=u_\theta+u_d\quad\text{in}\ \mathbf{H}^1_0(\Omega).
\]

The estimate \cref{energy_estimates} can now be obtained directly from this
decomposition.
By the coercivity of $\mathcal A_s$ and the definition of the
$\nabla_s$-operator,
\[
    \|\mathbf u_\theta\|_{\mathbf H_0^1(\Omega)}
    \lesssim \|\theta_s\|_{L^2(\Omega)}.
\]
Moreover, by the continuity of $\mathcal A_s^{-1}$ and $\mathcal B$,
\[
    \|\mathbf u_d\|_{\mathbf H_0^1(\Omega)}
    \lesssim
        \|f_s\|_{L^2(\Omega)}
        +
        \|\mathcal{S}_f\mathcal F_f\|_{\mathbf H_\Gamma^1(\Omega)}.
\]
The coercivity of $\mathcal A_f+\frac1{Da}I$ and the structure of
$\mathcal F_f$ yield
\[
    \|\mathcal{S}_f\mathcal F_f\|_{\mathbf H_\Gamma^1(\Omega)}
    \lesssim
    \|\mathcal F_f\|_{(\mathbf H_\Gamma^1(\Omega))^*}
    \lesssim
        \|f_f\|_{L^2(\Omega)}
        +
        \|f_N+\frac{\varphi_fa_2}{a_1}\mathbf n\|_{L^2(\Gamma_N)}.
\]
Since $a_1$, $a_2$, and $\varphi_f$ are fixed parameters and
$|\mathbf n|=1$ on $\Gamma_N$, we have
\[
\left\|
f_N+\frac{\varphi_fa_2}{a_1}\mathbf n
\right\|_{L^2(\Gamma_N)}
\lesssim
\|f_N\|_{L^2(\Gamma_N)}+1.
\]
Consequently,
\[
    \|\mathbf u_s\|_{\mathbf H_0^1(\Omega)}
    \lesssim
        \|\theta_s\|_{L^2(\Omega)}
        +
        \|f_s\|_{L^2(\Omega)}
        +
        \|f_f\|_{L^2(\Omega)}
        +
        \|f_N\|_{L^2(\Gamma_N)}+1,
\]
which proves \cref{energy_estimates}.

% The estimate \cref{energy_estimates} now can be established via standard energy estimates or alternatively by directly using the continuity of the involved operators: Using the coerciveness of $\mathcal{A}_s$ and the definition of the $\nabla$-operator, we have
% %
% \[
% c\|\mathbf{u}_s\|^2_{\mathbf{H}^1_0(\Omega)}\leq\langle\mathcal{A}_s\mathbf{u}_s,\mathbf{u}_s\rangle_{\mathbf{H}^1_0(\Omega)}=\int_\Omega \gamma_s\vartheta_s\nabla\cdot \mathbf{u}_s\di{x}+\int_\Omega\tilde{\mathcal{F}}\cdot \mathbf{u}_s\di{x}
% \]
% which leads to (using Young's inequality)
% %
% \[
% \|\mathbf{u}_s\|^2_{\mathbf{H}^1_0(\Omega)}\leq C\left(\|\theta_s\|^2_{L^2(\Omega)}+\|\tilde{\mathcal{F}}\|^2_{L^2(\Omega)}\right).
% \]
% %
% The coerciveness of $\mathcal{A}_f$ and the structure of $\mathcal{F}_f$ implies
% %
% \[
% c\|\mathcal{A}_f^{-1}\mathcal{F}_f\|^2_{\mathbf{H}^1_{\Gamma}(\Omega)}\leq \|\mathcal{F}_f\|^2_{(\mathbf{H}^1_{\Gamma}(\Omega))^*}\leq \|f_f\|^2_{L^2(\Omega)}+\|f_N+P_a\|^2_{L^2(\Gamma_N)}
% \]
% %
% leading to (using the continuity of $\mathcal{B}$)
% %
% \[
% \|\tilde{\mathcal{F}}\|^2_{L^2(\Omega)}\leq C\left(\|f_s\|^2_{L^2(\Omega)}+\|f_f\|^2_{L^2(\Omega)}+\|f_N+P_a\|^2_{L^2(\Gamma_N)}\right).
% \]
% %
% Finally, with \cref{eq:vf} the estimate \cref{energy_estimates} follows.
\end{proof}

With \Cref{lemma_step1}, we have the existence of a solution operator $\mathcal{L}\colon L^2(\Omega)\to \mathbf{H}^1_{\Gamma}(\Omega)\times \mathbf{H}^1_{0}(\Omega)$ mapping a given solid temperature profile $\theta_s$ to the corresponding hydrodynamic solution $(\mathbf{v}_f,\mathbf{u}_s)$.
In the fully-coupled problem (\cref{weak1,weak2,weak3,weak4,weak5}), the temperature profiles will depend on time as well and so, as a consequence, will the hydrodynamic solution.
Moreover, the regularity of $(\mathbf{v}_f,\mathbf{u}_s)$ with respect to time is essential for the full problem, specifically for $\mathbf{u}_s$ where we want to make sense of $\nabla\cdot\partial_t\mathbf{u}_s$ (\cref{weak5}).

In the following lemma, we collect some important regularity results for the solution $(\mathbf{v}_f,\mathbf{u}_s)$ with respect to a time parameter introduced via the data.
%%%%%%%%%%%%%%%%%%%%%%%%%%%%%%%%%%%%%%%%%%%%%%%%
% Lemma: Time regularity
%%%%%%%%%%%%%%%%%%%%%%%%%%%%%%%%%%%%%%%%%%%%%%%%
\begin{lemma}\label{lemma_time_regularity}
    For $t\in S$, let $t\mapsto(f_s(t),f_f(t),f_N(t),\theta_s(t))$ denote time-dependent data and let $t\mapsto(\mathbf{v}_f(t),\mathbf{u}_s(t))$ be the corresponding time-parametrized solution due to \cref{lemma_step1}.
    \begin{itemize}
        \item[$(i)$]  Let $f_s, f_f\in L^2(S\times\Omega)^d$, $\theta_s\in L^2(S\times\Omega)$, $f_N\in L^2(S\times\Gamma_N)^d$.
        Then,
        \[
        \mathbf{v}_f\in L^2(S;\mathbf{H}^1_{\Gamma}(\Omega)),\quad \mathbf{u}_s\in L^2(S;\mathbf{H}^1_{0}(\Omega)).
        \]
        \item[$(ii)$]If
        \[
            f_f\in L^\infty(S;L^2(\Omega)^d),
            \qquad
            f_N\in L^\infty(S;L^2(\Gamma_N)^d),
        \]
        then
        \[
            \mathbf v_f
            \in L^\infty(S;\mathbf H_\Gamma^1(\Omega))
            \hookrightarrow
            L^\infty(S;L^6(\Omega)^d).
        \]
        \item[$(iii)$] Let $\theta_s=0$. If $\partial_tf_f\in L^2(S;(\mathbf{H}_{\Gamma}^{1}(\Omega))^*)$, $\partial_tf_N\in L^2(S\times\Gamma_N)^d$, and $\partial_tf_s\in L^2(S;\mathbf{H}^{-1}(\Omega))$ we have 
        \[
        \partial_t\nabla\cdot \mathbf{u}_s\in L^2(S\times\Omega).
        \]
    \end{itemize}  
\end{lemma}
\begin{proof}
    $(i)$. The  $L^2$-time regularity of $(\mathbf{v}_f,\mathbf{u}_s)$ transfers in a straightforward way from the data as all involved operators are linear, continuous, and time-independent.
    
    $(ii)$. The estimate \cref{energy_estimates_a} applied for almost every
    $t\in S$ yields
    \[
        \|\mathbf v_f(t)\|_{\mathbf H_\Gamma^1(\Omega)}
        \lesssim
        \|f_f(t)\|_{L^2(\Omega)}
        +
        \|f_N(t)\|_{L^2(\Gamma_N)}
        +1.
    \]
    Hence,
    \[
        \mathbf v_f\in L^\infty(S;\mathbf H_\Gamma^1(\Omega)).
    \]
    Since $d\leq 3$, the Sobolev embedding
    $H^1(\Omega)\hookrightarrow L^6(\Omega)$ gives
    \[
        \mathbf v_f\in L^\infty(S;L^6(\Omega)^d).
    \]
    
    $(iii)$. With $\partial_tf_f\in L^2(S;(\mathbf{H}_{\Gamma}^{1}(\Omega))^*)$ and $\partial_tf_N\in L^2(S\times\Gamma_N)^d$ and taking into account the structure of $\mathcal{F}_f$ (\cref{operator5}), it holds
    \[
    \partial_t\mathcal{F}_f\in L^2(S;(\textbf{H}_\Gamma(\Omega))^*),\quad
    \langle\partial_t\mathcal{F}_f(t),w\rangle_{\textbf{H}_\Gamma(\Omega)}=\langle\partial_tf_f(t),w\rangle_{\textbf{H}_\Gamma(\Omega)}+\int_{\Gamma_N}\partial_tf_N(t)\cdot w\di{x}
    \]
    which implies, via \cref{lemma_step1}, $\partial_t\mathbf{v}_f\in L^2(S;\textbf{H}_\Gamma(\Omega))$ and, in consequence, $\partial_t\mathbf{u}_s\in L^2(S;\textbf{H}^1_0(\Omega))$.
\end{proof}

% \subsection{Analysis with solid drag $(\zeta>0)$}
% The problem takes the form 
% \begin{empheq}[left=(P)\empheqlbrace]{align*}
% a_f(\vf,w_f)+c_K(\vf-\zeta\partial_t\us,w_f)+b_f(p,w_f)&=l_f(t;w_f),\\
% a_s(\us,w_s)-c_K(\vf-\zeta\partial_t\us,w_s)+b_s(p,w_s)&=l_s(t;w_s),\\
% a_m(\vf,\partial_t\us,w)+m(p,q)&=l_m(t;q)
% \end{empheq}
% for all $(w_f,w_s,q)\in\Vf\times\Vs\times L^2(\Omega)$.
% %
% We approximate this problem by introducing $\us^i$ and $\partial_t\us(t_i)\approx\delta\us^i=\frac{\us(t_{i+1})-\us(t_i)}{\Delta t}$ and get
% %
% \begin{empheq}[left=(P_i)\empheqlbrace]{align*}
% a_f(\vf^i,w_f)+c_K(\vf^i-\zeta\delta\us^i,w_f)+b_f(p^i,w_f)&=l_f(t;w_f),\\
% a_s(\us^i,w_s)-c_K(\vf^i-\zeta\delta\us^i,w_s)+b_s(p^i,w_s)&=l_s(t;w_s),\\
% a_m(\vf^i,\delta\us^i,w)+m(p^i,q)&=l_m(t;q)
% \end{empheq}
% for all $(w_f,w_s,q)\in\Vf\times\Vs\times L^2(\Omega)$ and all $i=1,...,N$.

\subsection{Heat problem}
We now take a closer look at the thermo system whose weak formulation is given by \cref{weak4,weak5}.
In this problem, we have to account for the interphase heat coupling,
represented by the terms
\[
N(\theta_f-W\theta_s)
\quad\text{and}\quad
\kappa N\left(\theta_s-\frac{1}{W}\theta_f\right),
\]
as well as the thermo-mechanical coupling.
In the solid phase, it is given by $\delta'_s\nabla\cdot\partial_t\mathbf u_s$, and in the fluid phase, this consists of the convection term
$\nabla\cdot(\theta_f\mathbf v_f)$ and the effective divergence-coupling
term $\delta'_f\nabla\cdot\mathbf v_f$, where the latter contains both
the thermo-mechanical contribution and the contribution resulting from
the shift by the ambient temperature.

We first focus on the fluid heat system described by \cref{weak4}; this one is straightforward as we do not have to deal with the mixed derivative term (time and spatial) $\left(\nabla\cdot \mathbf{u}_s,\partial_t\varphi_s\right)_\Omega$ that comes up in the solid heat equation.
%
%%%%%%%%%%%%%%%%%%%%%%%%%%%%%%%%%%%%%%%%%%%%%%%%
% Lemma: Fluid heat problem
%%%%%%%%%%%%%%%%%%%%%%%%%%%%%%%%%%%%%%%%%%%%%%%%
\begin{lemma}\label{l:s2}
Let $g_f\in L^2(S\times\Omega)$ and $\theta_{f,0}\in L^2(\Omega)$.
For any $\theta_s\in L^2(S\times\Omega)$ and $ \mathbf v_f\in
    L^\infty(S;L^6(\Omega)^d)\cap L^2(S;\mathbf H_\Gamma^1(\Omega))$, there is a unique function $\theta_f\in \mathcal{W}$ satisfying
\begin{multline*}
-\int_S\int_\Omega
\mathrm{Pe}_f\theta_f\partial_t\psi_f
\di{x}\di{t}
+
\int_S\int_\Omega
\delta'_f\nabla\cdot\mathbf v_f\,\psi_f
\di{x}\di{t}
\\
+
\int_S\int_\Omega
\nabla\theta_f\cdot\nabla\psi_f
\di{x}\di{t}
-
\int_S\int_\Omega
\mathrm{Pe}_f\theta_f\mathbf v_f\cdot\nabla\psi_f
\di{x}\di{t}
\\
+
\int_S\int_\Omega
N(\theta_f-W\theta_s)\psi_f
\di{x}\di{t}
=
\int_S\int_\Omega
g_f\psi_f
\di{x}\di{t}
+
\int_\Omega
\mathrm{Pe}_f\theta_{f,0}\psi_f(0)
\di{x}
\end{multline*}
for all
\[
    \psi_f\in \left\{
        \psi\in L^2(S;H_0^1(\Omega))
        \cap H^1(S;L^2(\Omega))
        :
        \psi(T)=0
    \right\}.
\]
Moreover,
\begin{multline}
\|\theta_f\|^2_{L^\infty(S;L^2(\Omega))}+\|\partial_t\theta_f\|_{L^2(S;H^{-1}(\Omega))}^2+\|\nabla\theta_f\|^2_{L^2(S\times\Omega)}\\
\lesssim C_{\mathbf v_f}\left(\|\theta_s\|^2_{L^2(S\times\Omega)}+\|g_f\|^2_{L^2(S\times\Omega)}+\|\theta_{f,0}\|^2_{L^2(\Omega)}+
    \|\mathbf v_f\|_{L^2(S;\mathbf H_\Gamma^1(\Omega))}^2\right)
\end{multline}
where $C_{\mathbf v_f}>0$ depends continuously on
$\|\mathbf v_f\|_{L^\infty(S;L^6(\Omega))}$.
\end{lemma}
\begin{proof}
The convection term defines a continuous bilinear form on
$H_0^1(\Omega)$ since, for $d\leq3$,
\[
    \left|
    \int_\Omega
    \theta\,\mathbf v_f\cdot\nabla\psi\,\di{x}
    \right|
    \lesssim
    \|\mathbf v_f\|_{L^6(\Omega)}
    \|\theta\|_{H^1(\Omega)}
    \|\psi\|_{H^1(\Omega)}.
\]
Moreover, interpolation and Young's inequality yield
\[
    \left|
    \int_\Omega
    \theta\,\mathbf v_f\cdot\nabla\theta\,\di{x}
    \right|
    \leq
    \varepsilon\|\nabla\theta\|_{L^2(\Omega)}^2
    +
    C_\varepsilon
    \|\mathbf v_f\|_{L^6(\Omega)}^4
    \|\theta\|_{L^2(\Omega)}^2.
\]
Since
$\mathbf v_f\in L^\infty(S;L^6(\Omega)^d)$,
the standard theory for linear parabolic equations with
time-dependent bilinear forms yields existence and uniqueness, we refer to \cite[Chapter 7.1]{Evans2010}.
\end{proof}

Now, let $\mathbf{v}_f$ be the solution given by \cref{lemma_step1} and the data $f_f$, $f_N$ be as regular as required in \cref{lemma_time_regularity}.
Then, by \cref{lemma_time_regularity}.$(ii)$,
\[
    \mathbf v_f
    \in
    L^\infty(S;L^6(\Omega)^d)
    \cap
    L^2(S;\mathbf H_\Gamma^1(\Omega)),
\]
and hence \cref{l:s2} is applicable.
We introduce the solution operator $\mathcal{T}_{f}\colon L^2(S\times\Omega)\to \mathcal{W}$ assigning to a function $\theta_s\in L^2(S\times\Omega)$ the corresponding unique solution $\theta_f\in \mathcal{W}$ given by virtue of \Cref{l:s2} when $\mathbf{v}_f=\mathbf{v}_f$.
Note that $T_f-T_f(0)$ is linear, continuous, and, by Aubin-Lyons-Simon lemma \cite{Simon1986}, also compact from $L^2(S\times\Omega)$ to $L^2(S\times\Omega)$.
In that setting, any solution $(\theta_f,\theta_s)$ of the thermo-problem \cref{weak4,weak5} has to satisfy $\theta_f=\mathcal{T}_f\theta_s$.

Now, we start with the characterization of $\mathbf{u}_s$ given by \cref{U_operator} via \cref{lemma_step1} and splitting it in terms of a prescribed solid temperature $\theta_s$ and data parts, i.e.,
\[
\mathbf{u}_s(t)=u_\theta(t)+u_d(t) \quad\text{in}\ \mathbf{H}_0^1(\Omega)\ \text{for a.a.}\ t\in S.
\]
When $f_s,f_f, f_N$ satisfy Assumptions (A2) and (A3), we have $\partial_t\nabla\cdot U_d\in L^2(S\times\Omega)$ (\cref{lemma_time_regularity}(iii)) and can therefore slightly rearrange \cref{weak5} into
\begin{multline*}
-\int_S\int_\Omega
\mathrm{Pe}_s\theta_s\partial_t\psi_s
\di{x}\di{t}
+
\int_S\int_\Omega
\delta'_s
\nabla\cdot
\mathcal A_s^{-1}\nabla\theta_s\,
\partial_t\psi_s
\di{x}\di{t}
\\
+
\int_S\int_\Omega
\nabla\theta_s\cdot\nabla\psi_s
\di{x}\di{t}
+
\int_S\int_\Omega
\kappa N
\left(
    \theta_s-\frac{1}{W}T_f\theta_s
\right)
\psi_s
\di{x}\di{t}
\\
=
\int_S\int_\Omega
g_s\psi_s
\di{x}\di{t}
-\int_S\int_\Omega\delta'_s\,\partial_t\nabla\cdot\mathbf u_d\,\psi_s\di{x}\di{t}
\\
+
\int_\Omega
\left(
    \mathrm{Pe}_s\theta_{s,0}
    -
    \delta'_s
    \nabla\cdot
    \mathcal A_s^{-1}\nabla\theta_{s,0}
\right)
\psi_s(0)
\di{x}.
\end{multline*}
for all $\psi_s\in L^2(S;H_0^1(\Omega))$ with $\partial_t\psi_s\in L^2(S\times\Omega)$ and $\psi_s(T)=0$.
Introducing the linear operators
\begin{alignat*}{2}
\mathcal{K}&\colon H_0^1(\Omega)\to H^{-1}(\Omega),  &\qquad\langle\mathcal{K}\vartheta,\varphi\rangle_{H_0^1(\Omega)}&=\int_\Omega\nabla\vartheta\cdot\nabla\varphi\di{x},\\ 
\mathcal{R}&\colon L^2(\Omega)\to L^2(\Omega),&\quad \mathcal{R}&=-\nabla\cdot\mathcal{A}_s^{-1}\nabla_s,\\
\mathcal{M}&\colon L^2(\Omega)\to L^2(\Omega),&\quad \mathcal{M}&=
    \mathrm{Pe}_s\,\mathrm{Id}
    +
    \delta'_s\mathcal R.
\end{alignat*}
and the data $\mathcal{G}\in L^2(S\times\Omega)$ given by $\mathcal{G}=g_s-\delta'_s\partial_t\nabla\cdot\mathbf u_d$, the solid-heat system can be abstractly formulated as:
For $\theta_f\in L^2(S\times\Omega)$, find $\theta_s\in L^2(S;H^1_0(\Omega))$ with $\partial_t\left(\mathcal{M}\theta_s\right)\in L^2(S;H^{-1}(\Omega))$ and such that
\begin{equation}\label{operator_form}
    \partial_t(\mathcal M\theta_s)
    +
    \mathcal K\theta_s
    +
    \kappa N
    \left(
        \theta_s-\frac{1}{W}\theta_f
    \right)
    =
    \mathcal G
    \qquad
    \text{in }L^2(S;H^{-1}(\Omega)).
\end{equation}
Moreover, it has to satisfy $\lim_{t\to0}\mathcal{M}\theta_s(t)=\mathcal{M}\theta_{s,0}$ in $L^2(\Omega)$.

%%%%%%%%%%%%%%%%%%%%%%%%%%%%%%%%%%%%%%%%%%%%%%%%
% Lemma: Operator R
%%%%%%%%%%%%%%%%%%%%%%%%%%%%%%%%%%%%%%%%%%%%%%%%
\begin{lemma}\label{properties_of_R}
    The linear operator 
    \[
    \mathcal{R}\colon L^2(\Omega)\to L^2(\Omega),\quad \mathcal{R}=-\nabla\cdot\mathcal{A}_s^{-1}\nabla_s
    \]
    is continuous, self-adjoint, and positive.
\end{lemma}
\begin{proof}
Continuity follows from the continuity of the involved operators.
Let $f,g\in L^2(\Omega)$. For $\delta_s>0$,
\[
    \delta_s
    (\mathcal R f,g)_{L^2(\Omega)}
    =
    \left\langle
        \nabla g,
        \mathcal A_s^{-1}\nabla_s f
    \right\rangle.
\]
Since $\mathcal A_s^{-1}$ is symmetric, it follows that
\[
    (\mathcal R f,g)_{L^2(\Omega)}
    =
    (\mathcal R g,f)_{L^2(\Omega)}.
\]
Moreover,
\[
    \delta_s
    (\mathcal R f,f)_{L^2(\Omega)}
    =
    \left\langle
        \nabla f,
        \mathcal A_s^{-1}\nabla_s f
    \right\rangle
    \geq0,
\]
so that $\mathcal R$ is positive.
If $\delta_s=0$, then the operator $\nabla$ vanishes and hence $\mathcal R=0$.
\end{proof}
%%%%%%%%%%%%%%%%%%%%%%%%%%%%%%%%%%%%%%%%%%%%%%%%
% Lemma: Solid heat problem
%%%%%%%%%%%%%%%%%%%%%%%%%%%%%%%%%%%%%%%%%%%%%%%%
\begin{lemma}\label{l:s3}
Let $\theta_f\in L^2(S\times\Omega)$ as well as $\mathcal{G}\in L^2(S\times\Omega)$ and $\theta_{s,0}\in L^2(\Omega)$.
Then, there is a unique solution $\theta_s\in L^2(S;H^1_0(\Omega))$ with $\partial_t\left(\mathcal{M}\theta_s\right)\in L^2(S;H^{-1}(\Omega))$ satisfying $\lim_{t\to0}\mathcal{M}\theta_s(t)=\mathcal{M}\theta_{s,0}$ in $L^2(\Omega)$ and and solving \cref{operator_form}.
\end{lemma}
\begin{proof}
Observing that $\mathcal{R}$ is continuous, self-adjoint, and positive and $\mathcal{K}$ is coercive, the result for this implicit parabolic PDE follows from \cite[Chapter III, Proposition 3.2]{S96}.
\end{proof}
Moreover, since $\mathcal M$ is bounded, self-adjoint, and strictly
positive, the scalar product
\[
    (u,v)_{\mathcal M}
    :=
    (\mathcal Mu,v)_{L^2(\Omega)}
\]
is equivalent to the $L^2(\Omega)$ scalar product.
Hence, \cite[Chapter~III, Proposition~1.2]{S96}, together with
\[
    \theta_s\in L^2(S;H_0^1(\Omega)),
    \qquad
    \partial_t(\mathcal M\theta_s)
    \in L^2(S;H^{-1}(\Omega)),
\]
implies
\[
    \theta_s\in C([0,T];L^2(\Omega)).
\]
As a consequence, using the continuity of
\[
    \mathcal A_s^{-1}\nabla_s\colon
    L^2(\Omega)\to\mathbf H_0^1(\Omega),
\]
the temperature-induced displacement
\[
    \mathbf u_\theta
    =
    -\mathcal A_s^{-1}\nabla_s\theta_s
\]
belongs to $C([0,T];\mathbf H_0^1(\Omega))$.
In particular, $\mathbf u_\theta(0)\in \mathbf H_0^1(\Omega)$ is well defined.
With these results, we can now prove the general well-posedness:

\begin{theorem}[Existence result]\label{theorem:existence}
    Let Assumptions $(A1)$--$(A4)$ be satisfied.
    Then there exists a unique weak solution
    \[
        \mathbf v_f
        \in L^\infty(S;\mathbf H_\Gamma^1(\Omega)),
        \qquad
        \mathbf u_s
        \in C([0,T];\mathbf H_0^1(\Omega)),\qquad
        P\in L^2(S\times\Omega),
        \qquad
        \theta_f\in W(S,\Omega),
    \]
    and
    \[
        \theta_s\in L^2(S;H_0^1(\Omega)),
        \qquad
        \partial_t(\mathcal M\theta_s)
        \in L^2(S;H^{-1}(\Omega)).
    \]
    solving the non-dimensional problem in the sense of \cref{weak1,weak2,weak3,weak4,weak5}.
\end{theorem}
\begin{proof}
Set
\[
    X:=L^2(S\times\Omega).
\]
By \cref{l:s2}, for every prescribed $\theta_s\in X$ there exists a unique
fluid temperature $\theta_f=\mathcal T_f\theta_s$ where
\[
    \mathcal T_f:X\to X
\]
is an affine solution operator.
We write
\[
    \mathcal T_f\eta
    =
    \theta_f^0+\mathcal L_f\eta,
    \qquad
    \theta_f^0:=\mathcal T_f0.
\]
The operator $\mathcal L_f:X\to X$ is linear and continuous, and compact (which follows by the estimate in \cref{l:s2} together with the Aubin--Lions--Simon compactness theorem).
Similarly, \cref{l:s3} defines an affine solution operator
\[
    \mathcal T_s:X\to X,
    \qquad
    \mathcal T_s\eta
    =
    \theta_s^0+\mathcal L_s\eta,
    \qquad
    \theta_s^0:=\mathcal T_s0.
\]
The linear part $\mathcal L_s:X\to X$ is continuous. Indeed, for
$z=\mathcal L_s\eta$ we have
\[
    \partial_t(\mathcal Mz)
    +\mathcal Kz
    +\kappa Nz
    =
    \frac{\kappa N}{W}\eta,
    \qquad
    \mathcal Mz(0)=0,
\]
and the standard energy estimate yields
\[
    \|z\|_X\lesssim\|\eta\|_X.
\]
A solution of the coupled thermal problem therefore must satisfy
\[
    \theta_s
    =
    \mathcal T_s(\mathcal T_f\theta_s),
\]
or equivalently
\[
    (I-\mathcal C)\theta_s=h,
\]
where
\[
    \mathcal C:=\mathcal L_s\mathcal L_f,
    \qquad
    h:=\theta_s^0+\mathcal L_s\theta_f^0.
\]
Since $\mathcal L_f$ is compact and $\mathcal L_s$ is continuous,
$\mathcal C:X\to X$ is compact.
It remains to show that
\[
    \ker(I-\mathcal C)=\{0\}.
\]
To that end, let $\vartheta_s\in\ker(I-\mathcal C)$ and set
\[
    \vartheta_f:=\mathcal L_f\vartheta_s.
\]
Then $(\vartheta_f,\vartheta_s)$ solves the homogeneous coupled thermal
problem with zero initial data.
Testing the fluid equation with $\vartheta_f$ and the solid equation with $\frac{W^2}{\kappa}\vartheta_s$, and adding the resulting identities,
gives
\begin{multline*}
    \frac{\di}{\di t}
    \bigg(
        \frac{{\rm Pe}_f}{2}\|\vartheta_f\|_{L^2(\Omega)}^2
        +
        \frac{W^2}{2\kappa}
        (\mathcal M\vartheta_s,\vartheta_s)_{L^2(\Omega)}
    \bigg)+
    \|\nabla\vartheta_f\|_{L^2(\Omega)}^2
    +
    \frac{W^2}{\kappa}
    \|\nabla\vartheta_s\|_{L^2(\Omega)}^2
    \\
    +N\|\vartheta_f-W\vartheta_s\|_{L^2(\Omega)}^2
    =
    {\rm Pe}_f
    \int_\Omega
    \vartheta_f\mathbf v_f\cdot\nabla\vartheta_f
    \di{x}.
\end{multline*}
Here, for the solid contribution, we used the generalized product rule
for implicit evolution equations associated with the bounded,
self-adjoint, and positive operator $\mathcal M$, namely
\[
    \left\langle
        \partial_t(\mathcal M\vartheta_s),
        \vartheta_s
    \right\rangle
    =
    \frac{1}{2}\frac{\di}{\di t}
    (\mathcal M\vartheta_s,\vartheta_s)_{L^2(\Omega)},
\]
see \cite[Chapter~III, Proposition~1.2]{S96}.
Using
$\mathbf v_f\in L^\infty(S;L^6(\Omega)^d)$,
interpolation, and Young's inequality, we obtain
\[
    \left|
    \int_\Omega
    \vartheta_f\mathbf v_f\cdot\nabla\vartheta_f
    \di{x}
    \right|
    \leq
    \frac12\|\nabla\vartheta_f\|_{L^2(\Omega)}^2
    +
    C\|\mathbf v_f\|_{L^6(\Omega)}^4
    \|\vartheta_f\|_{L^2(\Omega)}^2.
\]
Moreover, since $\mathcal M$ is positive and
\[
    (\mathcal M z,z)_{L^2(\Omega)}
    \geq
    {\rm Pe}_s\|z\|_{L^2(\Omega)}^2,
\]
Gronwall's inequality and the homogeneous initial data imply
\[
    \vartheta_f=\vartheta_s=0.
\]
Thus,
\[
    \ker(I-\mathcal C)=\{0\}.
\]
The Fredholm alternative now implies that $I-\mathcal C$ is invertible
on $X$.
Hence, there exists a unique $\theta_s\in X$ satisfying the coupled thermal problem, and
\[
    \theta_f=\mathcal T_f\theta_s
\]
is uniquely determined. The additional regularity of the temperatures
follows from \cref{l:s3,l:s2}.

Finally, for almost every $t\in S$, Lemma~1 yields the corresponding
unique hydrodynamic variables $(\mathbf v_f,\mathbf u_s)$, and the
pressure is recovered from
\[
    P=\frac{1}{a_1}
    \left(
        a_2-\varphi_f\nabla\cdot\mathbf v_f
    \right).
\]
The asserted regularity follows from \cref{l:s2,lemma_time_regularity,lemma_step1} and the preceding
estimates. Uniqueness of the full solution follows from the uniqueness
of the thermal and hydrodynamic subproblems.
\end{proof}

\section{Reduced model}\label{sec:exact_solutions}

So far, we have proposed a coupled model of thermo-elasto-hydrodynamics and established the existence and uniqueness of a corresponding weak solution. 
In this section, we consider a simplified one-dimensional setting in order to further investigate the interplay between fluid flow, solid deformation, and heat transfer.
More precisely, we consider the non-dimensional system \cref{main-5a,main-5b,main-5c,main-5g,main-5h} in one spatial dimension and assume a stationary regime.
\begin{figure}[H]
\centering
\includegraphics[width=\textwidth]{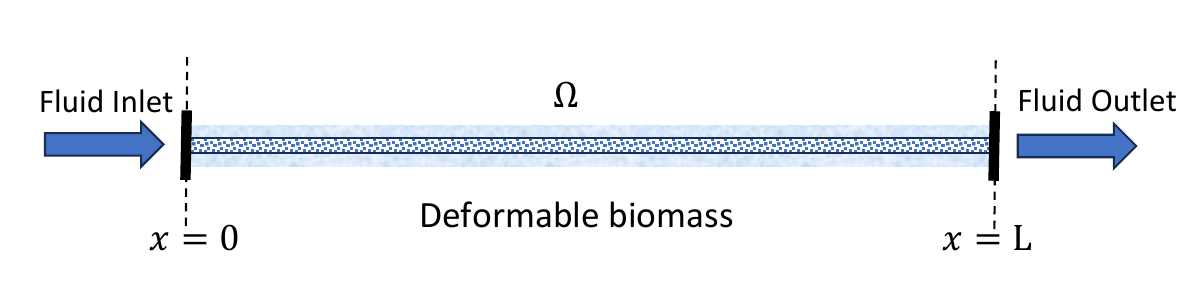}
\caption{Schematic representation of the 1D model.}\label{Figure1}
\end{figure}
Hence, all unknowns are assumed to be time-independent.
Moreover, the body-force densities and heat-production terms are neglected.
We denote the corresponding one-dimensional profiles by
\[
    V_f=V_f(x),\qquad
    U_s=U_s(x),\qquad
    P=P(x),\qquad
    \theta_f=\theta_f(x),\qquad
    \theta_s=\theta_s(x).
\]
For comparison purposes, we additionally introduce a parameter $\Xi\in\{0,1\}$ multiplying the thermoelastic contribution in the solid momentum balance. 
The case $\Xi=1$ corresponds to the original 
thermoelastic coupling, whereas $\Xi=0$ suppresses this contribution.
Under these simplifying assumptions, the one-dimensional system reads
\begin{subequations}
\begin{align}
    -(2+\Lambda_f)\frac{\di^2 V_f}{\di x^2}
    +\varphi_f\frac{\di P}{\di x}
    +\frac{1}{Da}V_f
    &=0,
    \label{F1}
    \\
    -(2+\Lambda_s)\frac{\di^2 U_s}{\di x^2}
    +\varphi_s\frac{\di P}{\di x}
    -\frac{1}{Da}V_f
    +\delta_s\Xi\frac{\di\theta_s}{\di x}
    &=0,
    \label{S1}
    \\
    \frac{\di}{\di x}(\varphi_f V_f)
    &=a_2-a_1P,
    \label{P1}
    \\
    \mathrm{Pe}_f
    \frac{\di}{\di x}(\theta_f V_f)
    -\frac{\di^2\theta_f}{\di x^2}
    &=
    -N(\theta_f-W\theta_s)
    -\delta_f'\frac{\di V_f}{\di x},
    \label{T1}
    \\
    -\frac{\di^2\theta_s}{\di x^2}
    &=
    -\kappa N
    \left(
        \theta_s-\frac{1}{W}\theta_f
    \right).
    \label{T2}
\end{align}
Here, we impose boundary conditions adapted to the inlet--outlet configuration shown in \cref{Figure1}.
These differ from the boundary conditions used in the well-posedness analysis of \cref{sec:analysis} and are chosen here for the subsequent semi-analytical study.

At the inlet, the fluid velocity is prescribed by the characteristic
velocity used in the non-dimensionalization, while at the outlet we
assume a vanishing velocity gradient. Hence,
\begin{align}
    V_f(0) = 1,\qquad
    \frac{\di V_f}{\di x}(1) = 0.
    \label{FB1}
\end{align}
For the solid displacement, we assume that the solid structure is fixed
at both ends,
\begin{align}
    U_s(0)=U_s(1) &= 0.
    \label{SB1}
\end{align}
For the thermal subsystem, we prescribe the temperatures at the inlet
and assume zero heat flux at the outlet. Thus,
\begin{align}
    \theta_s(0)=\theta_f(0) &= 1,\label{TBC1}
    \\
    \frac{\di\theta_s}{\di x}(1)=
    \frac{\di\theta_f}{\di x}(1) &= 0.
    \label{TBC2}
\end{align}
\end{subequations}
Using \cref{T2}, the fluid temperature can be expressed in terms of
the solid temperature as
\begin{align}
    \theta_f
    =
    W\theta_s
    -
    \frac{W}{\kappa N}
    \frac{\di^2\theta_s}{\di x^2}.
    \label{theta_f_from_s}
\end{align}
Substituting this relation into \cref{T1}, we obtain the fourth-order
equation
\begin{multline}
    \frac{\di^4\theta_s}{\di x^4}
    -
    \mathrm{Pe}_f V_f
    \frac{\di^3\theta_s}{\di x^3}
    -
    \left(
        (\kappa+1)N
        +
        \mathrm{Pe}_f\frac{\di V_f}{\di x}
    \right)
    \frac{\di^2\theta_s}{\di x^2}\\
    +
    \kappa N\mathrm{Pe}_f V_f
    \frac{\di\theta_s}{\di x}
    +
    \kappa N\mathrm{Pe}_f
    \frac{\di V_f}{\di x}\theta_s
    =
    -
    \frac{\kappa N\delta_f'}{W}
    \frac{\di V_f}{\di x}.
    \label{TTT1}
\end{multline}
The corresponding boundary conditions follow from
\cref{TBC1,TBC2} and \cref{theta_f_from_s} and are given by
\begin{subequations}
\begin{align}
    \theta_s(0) &= 1,
    \label{TB1}\\
    W\theta_s(0)
    -
    \frac{W}{\kappa N}
    \frac{\di^2\theta_s}{\di x^2}(0)
    &=1,
    \label{TB2}\\
    \frac{\di\theta_s}{\di x}(1) &=0,
    \label{TB3}\\
    W\frac{\di\theta_s}{\di x}(1)
    -
    \frac{W}{\kappa N}
    \frac{\di^3\theta_s}{\di x^3}(1)
    &=0.
    \label{TB4}
\end{align}
\end{subequations}

\paragraph{Solution procedure.}
We first consider the elastohydrodynamic subsystem \cref{F1,P1}.
From \cref{P1}, the pressure is given by
\begin{align}
    P
    =
    \frac{1}{a_1}
    \left(
        a_2
        -
        \varphi_f\frac{\di V_f}{\di x}
    \right).
    \label{pressure_1d}
\end{align}
Substituting this relation into \cref{F1} gives
\begin{align}
    \left(
        2+\Lambda_f+\frac{\varphi_f^2}{a_1}
    \right)
    \frac{\di^2V_f}{\di x^2}
    -
    \frac{1}{Da}V_f
    =0.
\end{align}
Introducing
\begin{align}
    \alpha_f^2
    :=
    \frac{1}{
        Da\left(
            2+\Lambda_f+\frac{\varphi_f^2}{a_1}
        \right)
    },
    \label{alpha_f}
\end{align}
the fluid velocity satisfies the ODE
\[
    \frac{\di^2V_f}{\di x^2}
    -
    \alpha_f^2V_f
    =0,\qquad V_f(0) = 1,\qquad
    \frac{\di V_f}{\di x}(1) = 0.
\]
It is easy to check that the solution is given via
\begin{align}
    V_f(x)
    =
    \frac{
        \cosh\bigl(\alpha_f(1-x)\bigr)
    }{
        \cosh(\alpha_f)
    }.
    \label{FS}
\end{align}
which directly implies
\begin{align}
    P(x)
    =
    \frac{1}{a_1}
    \left[
        a_2
        +
        \varphi_f\alpha_f
        \frac{
            \sinh\bigl(\alpha_f(1-x)\bigr)
        }{
            \cosh(\alpha_f)
        }
    \right].
    \label{PS}
\end{align}

Having determined the fluid velocity and pressure explicitly, we next consider the thermal subsystem.
In the stationary setting considered here, \cref{T1,T2} are independent of the solid displacement and can therefore be solved independently of the solid momentum balance.
With the fluid velocity now explicitly known, the fourth-order boundary-value problem \cref{TTT1} together with the corresponding boundary conditions is completely determined. 
In general, due to the spatially varying coefficients $V_f$ and $\di V_f/\di x$, we solve this problem numerically.
Once $\theta_s$ is obtained, the fluid temperature $\theta_f$ is recovered directly from \cref{theta_f_from_s}.

Finally, we consider the solid displacement.
From \cref{pressure_1d,alpha_f}, we have
\[
    \frac{\di P}{\di x}
    =
    -\frac{\varphi_f}{a_1}
    \frac{\di^2V_f}{\di x^2}
    =
    -\frac{\varphi_f\alpha_f^2}{a_1}V_f.
\]
Hence, \cref{S1} can be written as
\begin{align}
    \frac{\di^2U_s}{\di x^2}
    =
    \alpha_s V_f
    +
    \frac{\delta_s\Xi}{2+\Lambda_s}
    \frac{\di\theta_s}{\di x},
    \label{solid_reduced}
\end{align}
where
\begin{align}
    \alpha_s
    :=
    -\frac{1}{2+\Lambda_s}
    \left(
        \frac{\varphi_s\varphi_f\alpha_f^2}{a_1}
        +
        \frac{1}{Da}
    \right).
    \label{alpha_s}
\end{align}
The solution is given by
\begin{align}
    U_s(x)
    =
    \frac{\alpha_s}{\alpha_f^2}
    \left[
        V_f(x)-1
        +
        x\left(
            1-\frac{1}{\cosh(\alpha_f)}
        \right)
    \right]
    +
    \frac{\Xi\delta_s}{2+\Lambda_s}
    \left[
        \int_0^x\theta_s(y)\,\di y
        -
        x\int_0^1\theta_s(y)\,\di y
    \right]
    \label{US1}
\end{align}
where we note that $U_s$ does not depend on $\theta_s$ when $\Xi=0$.

To summarize, $V_f$ and $P$ are given explicitly by \cref{FS,PS}.
The solid temperature $\theta_s$ is determined numerically from the fourth-order boundary-value problem \cref{TTT1,TB1,TB2,TB3,TB4}, after which $\theta_f$ is recovered from \cref{theta_f_from_s} and $U_s$ from \cref{US1}.

\section{Results and Discussion}\label{sec:results_discussion}
We illustrate the behavior of the reduced one-dimensional model presented in \cref{sec:exact_solutionsb} and investigate the influence of selected model parameters on the fluid velocity, pressure, solid displacement, and temperature profiles.
Unless stated otherwise, the non-dimensional reference parameters used in the simulations are listed in \cref{t:parameters}.
To relate the porosity of the tissue to its permeability, we use a Kozeny--Carman-type relation.
Following \cite{majumder2023non}, the intrinsic permeability $k=\mu_f K$ is given by
\[
    \mu_f K
    =
    \frac{\phi_f^3}
         {C_k(1-\phi_f)^2 S^2},
\]
where $S$ denotes the wetted surface area per unit volume and $C_k$
is the Kozeny--Carman constant. Assuming cells of characteristic
diameter $D_c$ and using $S=4/D_c$, this becomes
\[
    \mu_f K
    =
    \frac{\phi_f^3D_c^2}
         {16C_k(1-\phi_f)^2}.
\]
\begin{table}[H]
\footnotesize
\centering
\begin{tabular}{llll}
\tline
Variable & Value & Description \\ \tline
$\phi_f$ & $0.9$ & Fluid volume fraction \\
$\phi_s$ & $0.1$ & Solid volume fraction \\
$Da$ & $0.0152$ & Darcy number \\
$a_1$ & $1$ & Pressure-source parameter \\
$a_2$ & $1.5$ & Pressure-source parameter \\
$\Lambda_f$ & $1$ & Fluid viscosity ratio \\
$\Lambda_s$ & $1$ & Solid elasticity ratio \\
$Pe_f$ & $0.7596$ & Fluid Péclet number \\
$N$ & $2$ & Interphase heat exchange parameter \\
$\kappa$ & $1$ & Thermal conductivity ratio \\
$W$ & $1$ & Temperature-scale ratio \\
$\delta_f'$ & $1$ & Fluid divergence-coupling parameter \\
$\delta_s$ & $0.1$ & Solid thermoelastic coupling parameter \\ \tline
\end{tabular}
\caption{Reference set of non-dimensional parameters used in the simulations.
Individual parameters are varied as specified in the corresponding figures.}
\label{t:parameters}
\end{table}
For the reference high-porosity case listed in \cref{t:parameters},
we use $D_c=10\,\mu\mathrm{m}$, $L=100\,\mu\mathrm{m}$, and
$C_k=3$, which yields $Da=0.0152$ for $\phi_f=0.9$.

\begin{figure}[h]
\centering
\includegraphics[width=.75\textwidth]{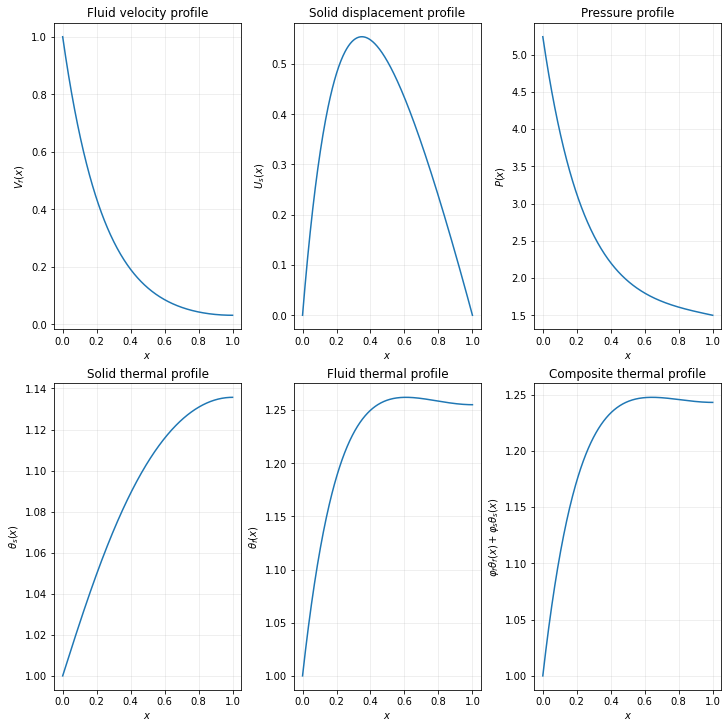}
\caption{For high porosity $\phi_f=0.9$, permeability $Da=0.0152$, and $\kappa=1$, $N=2$, $\mbox{Pe}_f=0.7596$.}\label{A2}
\end{figure}

\begin{figure}[h]
\centering
\includegraphics[width=0.75\textwidth]{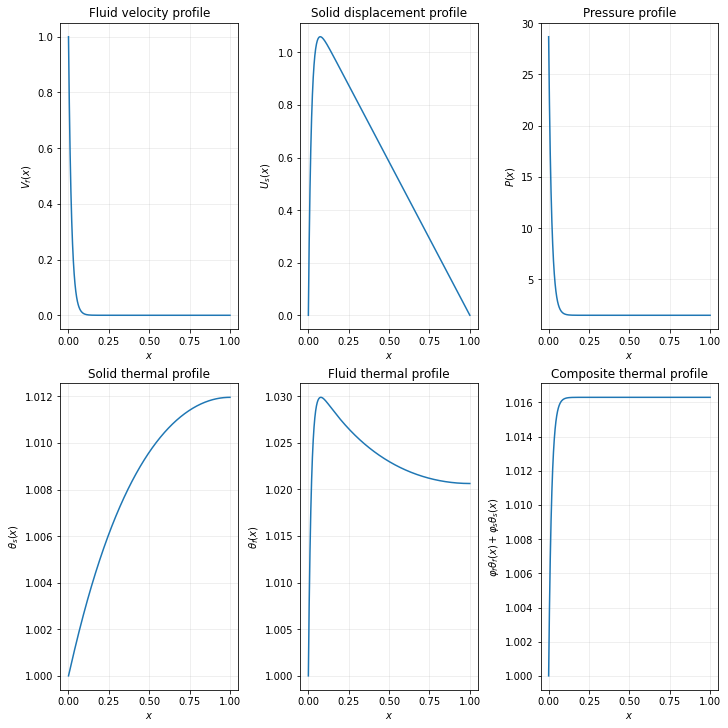}
\caption{For low porosity $\phi_f=0.5$ permeability,  $Da=0.000104$, and $\kappa=1$, $N=2$, $\mbox{Pe}_f=0.7596$.}\label{A3}
\end{figure}

In \cref{A2,A3}, we compare two porosity--permeability regimes linked through the Kozeny--Carman relation described above.
For $\phi_f=0.9$ we obtain $Da=0.0152$, whereas $\phi_f=0.5$ leads to the much smaller value $Da=0.000104$.
The lower permeability leads to a faster decay of the fluid velocity close to the inlet and, correspondingly, to a substantially larger pressure gradient.
The solid displacement is also more pronounced and its maximum is shifted towards the inlet due to the more localized hydrodynamic forcing.
A similar effect is also visible in the thermal profiles:
In the high-permeability case, the fluid temperature varies over a substantial part of the domain and the solid temperature follows through interphase heat exchange.
For the low-permeability case, the fluid-temperature variation is concentrated near the inlet and the solid temperature remains comparatively weakly varying.

\begin{figure}[h]
\centering
\includegraphics[width=.75\textwidth]{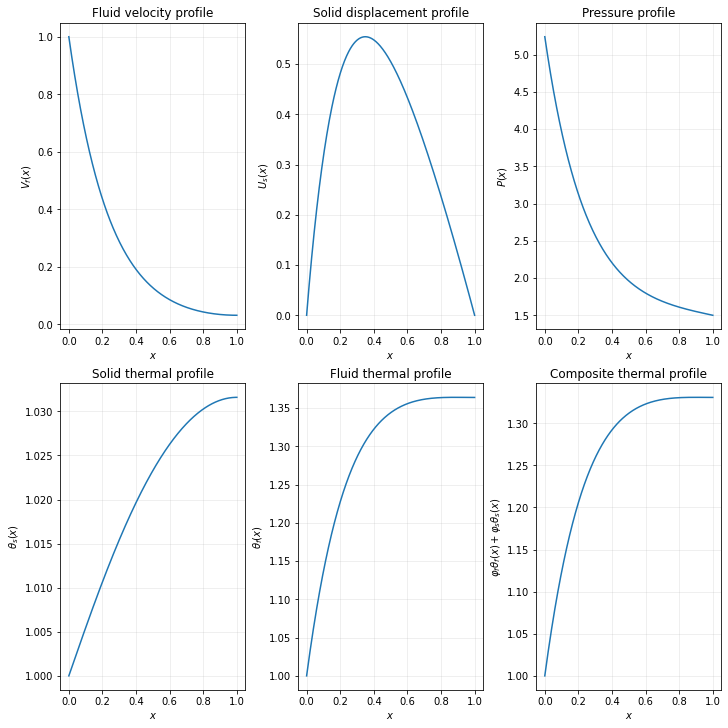}
\caption{Weak thermic coupling $N=0.2$ and $\phi_f=0.9$, $Da=0.0152$, $\kappa=1$, $\mbox{Pe}_f=0.7596$.}\label{A4}
\end{figure}

\begin{figure}[h]
\centering
\includegraphics[width=.75\textwidth]{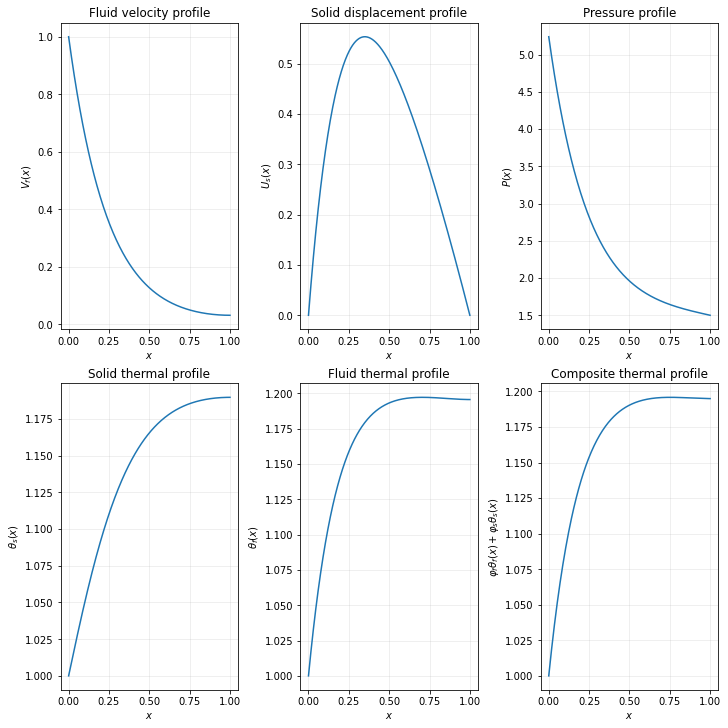}
\caption{Strong thermic coupling $N=20$ and $\phi_f=0.9$, $Da=0.0152$, $\kappa=1$, $\mbox{Pe}_f=0.7596$.}\label{A5}
\end{figure}

Next, we investigate the influence of the interphase heat exchange
parameter $N$, while keeping the porosity--permeability regime and all remaining parameters fixed (see \cref{A4,A5}). 
Comparing the weakly coupled case $N=0.2$ with the strongly coupled case $N=20$, we observe that the hydrodynamic profiles remain unchanged.
The main effect of increasing $N$ is a substantially stronger thermal coupling between the two phases.
For small $N$, the fluid temperature rises considerably above the solid temperature, whereas for large $N$ the two temperature profiles remain much closer to each other.
As expected, the stronger heat exchange reduces the temperature difference between the phases.

\begin{figure}[h]
\centering
\includegraphics[width=.75\textwidth]{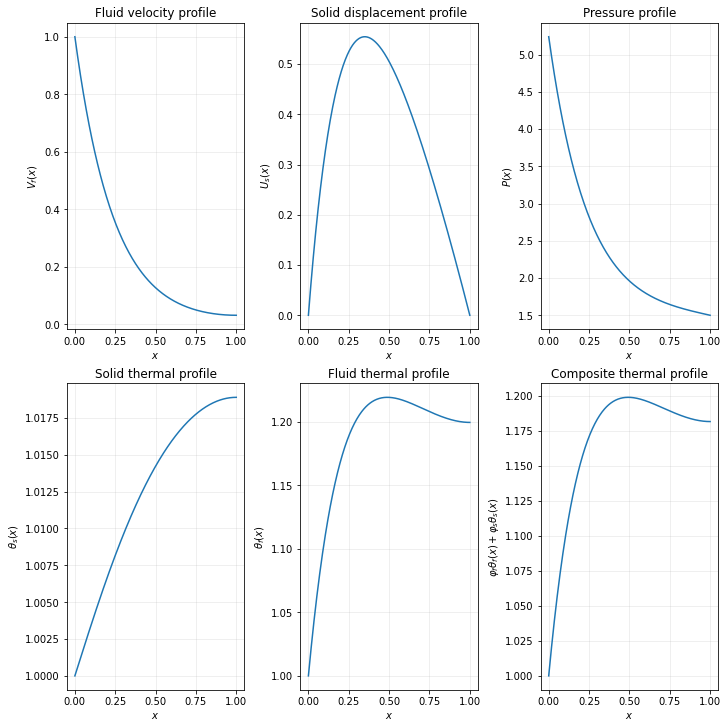}
\caption{Low conductivity ratio: $\kappa=0.1$ and $\phi_f=0.9$, $Da=0.0152$, $N=2$, $\mbox{Pe}_f=0.7596$.}\label{A6}
\end{figure}

\begin{figure}[h]
\centering
\includegraphics[width=.75\textwidth]{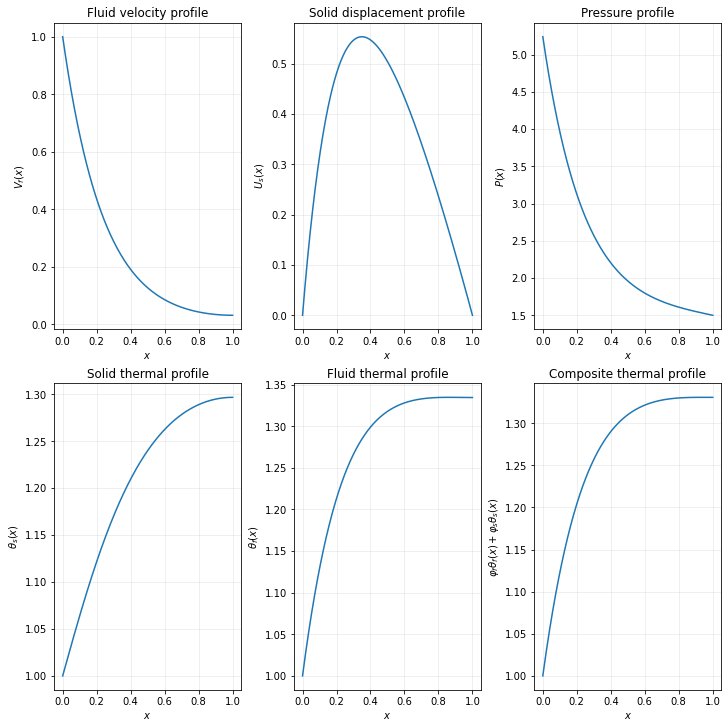}
\caption{High conductivity ratio: $\kappa=5$ and $\phi_f=0.9$, $Da=0.0152$, $N=2$, $\mbox{Pe}_f=0.7596$.}\label{A7}
\end{figure}

Finally, we investigate the influence of the thermal conductivity ratio $\kappa=\kappa_f/\kappa_s$ while keeping the remaining parameters fixed (see \cref{A6,A7}).
For small $\kappa$, the solid temperature responds only weakly to the fluid temperature and a pronounced temperature difference between the two phases remains. 
Increasing $\kappa$ strengthens the thermal response of the solid phase, so that $\theta_s$ follows $\theta_f$ more closely.
As expected, the hydrodynamic profiles are unaffected by this variation.
These examples illustrate how the hydrodynamic and thermal coupling
mechanisms can be separated and investigated within the reduced model.

\FloatBarrier

\section*{Acknowledgments}
The research activity of ME was partially funded by the European Union’s Horizon 2022 research and innovation program under the Marie Sklodowska-Curie fellowship project MATT (funding number 101061956, \url{https://doi.org/10.3030/101061956}).

\section*{Conflict of interest}
On behalf of all authors, the corresponding author states that there is no conflict of interest. 

\section*{Data availability statement}
Data will be made available on request.

\bibliographystyle{abbrv}
\bibliography{biblio_porous_medium}

\end{document}